\documentclass{CustomArticle}

\usepackage{FramedTheoremsForCustomArticle}
\usepackage{CommonMath}
\usepackage{ProjectMacros}

\usepackage{algorithm}
\usepackage{algorithmic}
\usepackage[
  style=trad-abbrv,
  citestyle=numeric-comp,
  sorting=nyt
]{biblatex}
\usepackage{booktabs}
\usepackage[font=small]{caption}
\usepackage{graphicx}
\usepackage{import}
\usepackage{microtype}

\usepackage{hyperref}
\usepackage[capitalize,compress]{cleveref}

\PaperGeneralInfo{
  title = {\vspace*{-2cm} {Polynomial Preconditioning for Gradient Methods}},
  date = {January 30, 2023}
}
\AddPaperAuthor{
  name = {Nikita Doikov},
  affiliation = {EPFL, Switzerland},
  email = {nikita.doikov@epfl.ch}
}
\AddPaperAuthor{
  name = {Anton Rodomanov},
  affiliation = {Catholic University of Louvain (UCLouvain), Belgium},
  email = {anton.rodomanov@uclouvain.be}
}
\PaperAbstract{%
  We study first-order methods with preconditioning for solving structured
  nonlinear convex optimization problems.
  We propose a new family of preconditioners generated by symmetric polynomials.
  They provide first-order optimization methods with a provable improvement of
  the condition number, cutting the gaps between highest eigenvalues, without
  explicit knowledge of the actual spectrum.
  We give a stochastic interpretation of this preconditioning in terms of
  coordinate volume sampling and compare it with other classical approaches,
  including the Chebyshev polynomials.
  We show how to incorporate a polynomial preconditioning into the Gradient and
  Fast Gradient Methods and establish the corresponding global complexity
  bounds.
  Finally, we propose a simple adaptive search procedure that automatically
  chooses the best possible polynomial preconditioning for the Gradient Method,
  minimizing the objective along a low-dimensional Krylov subspace.
  Numerical experiments confirm the efficiency of our preconditioning strategies
  for solving various machine learning problems.
}
\PaperKeywords{%
  Convex optimization,
  Preconditioning,
  Gradient methods,
  Accelerated methods,
  Symmetric polynomials
}
\PaperThanks{%
  The work of the first author was supported by the Swiss State Secretariat for
  Education, Research and Innovation (SERI) under contract number~22.00133.
  The work of the second author received funding from the European Research
  Council (ERC) under the European Union's Horizon 2020 research and innovation
  programme (grant agreement No.~788368)
}

\begin{document}

  \PrintTitleAndAbstract

  \section{Introduction}
\label{SectionIntro}

\paragraph{Motivation. }

Preconditioning is an important tool for improving the performance of
numerical algorithms. The classical example
is the preconditioned \textit{Conjugate Gradient Method}
\cite{hestenes1952methods} for solving a system of linear equations.
It proposes to modify the initial system in a way
to improve its eigenvalue distribution and thus to accelerate the convergence of the method.
The question of choosing the right preconditioner heavily depends on the problem structure,
and there exist many problem-specific recommendations
which provide us with a good trade-off between computational cost
and the spectrum properties of the new system. Some notable examples include
\textit{Jacobi} or the \textit{diagonal} preconditioners,
\textit{symmetric successive over-relaxation},
the \textit{incomplete Cholesky} factorization \cite{golub2013matrix},
\textit{Laplacian} preconditioning for graph problems
\cite{vaidya1991solving,spielman2004nearly},
preconditioners for
discretizations of system
of \textit{partial differential equations} \cite{mardal2011preconditioning}.

Another important class of numerical algorithms are the second-order methods or
\textit{Newton's Method} (see, e.g. \cite{nesterov2018lectures}),
that aims to solve difficult ill-conditioned problems by
using local curvature information (the Hessian matrix)
as a preconditioner at every step.
However, being a powerful algorithm, each iteration of Newton's Method
is very expensive. It requires to solve a system of linear equations with the Hessian matrix,
and in case of quadratic objective it is equivalent to solving the original problem.

\begin{figure}[h]
	\centering
	\includegraphics[width=0.70\linewidth]{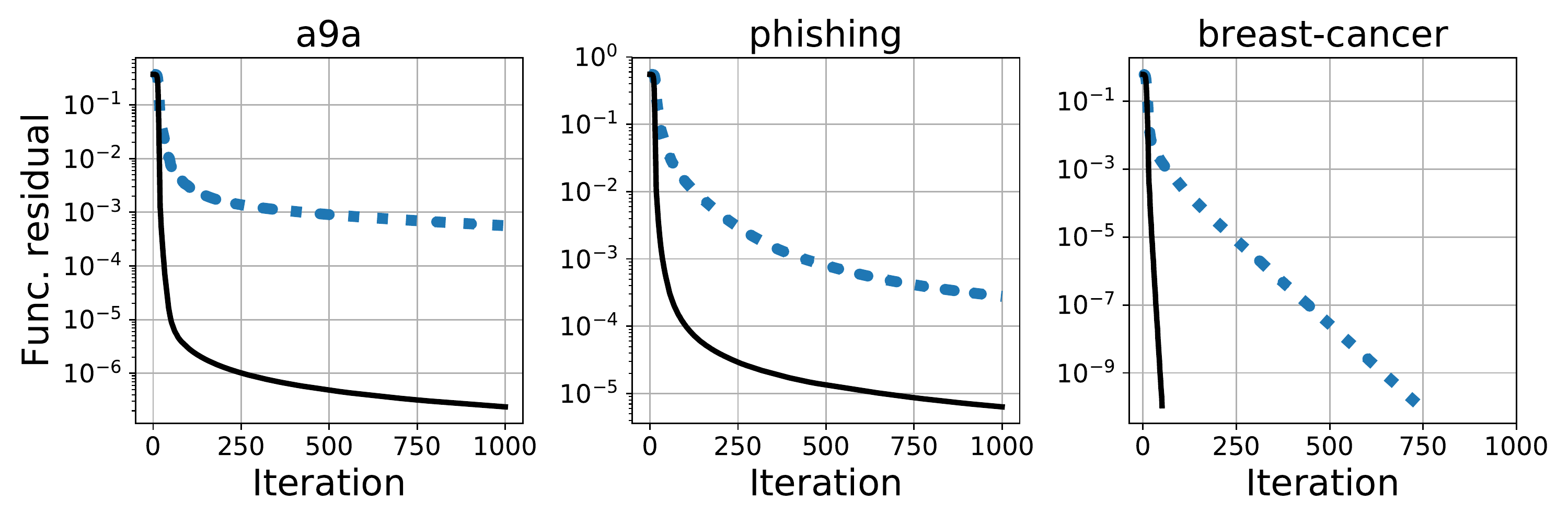}
	\includegraphics[width=0.70\linewidth]{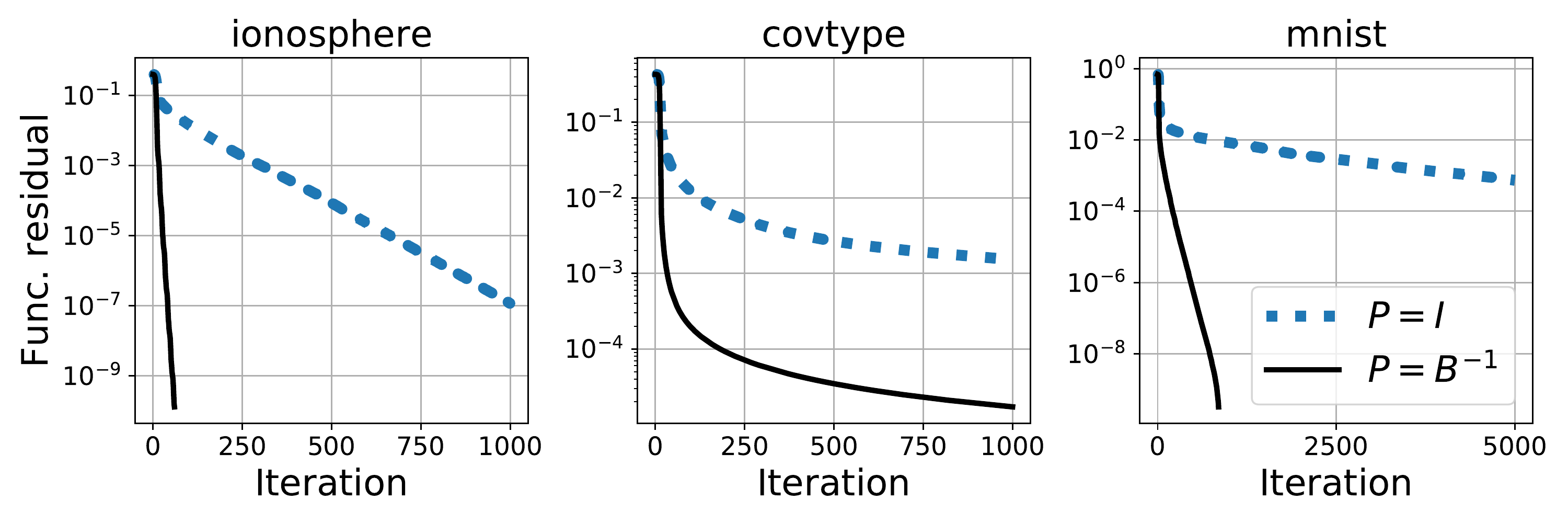}

	\caption{\small Training  logistic regression with the standard gradient descent ($\mat{P} = \mat{I}$),
		and using the inverse of the curvature matrix $(\mat{P} = \mat{B}^{-1})$
		as a preconditioner in \eqref{IntroGM}.
		The latter method works much faster, while it can be very expensive
		to compute $\mat{B}^{-1}$ for large scale problems.
	}
	\label{fig:Binv}
\end{figure}

In this paper, our goal is to solve a general \textit{nonlinear optimization problem}
with a structured convex objective by the efficient first-order methods.
Thus, in the case of unconstrained minimization
of a smooth function: $\min_{\xx} f(\xx)$, the simplest method that we study
is as follows,
for $k \geq 0$:
\beq \label{IntroGM}
\ba{rcl}
\xx_{k + 1} & = & \xx_k - \alpha_k \mat{P} \nabla f(\xx_k),
\ea
\eeq
where $\alpha_k > 0$ is a stepsize and $\mat{P}$ is a fixed preconditioning matrix.
$\mat{P} := \mat{I}$ corresponds to the classical gradient descent.
Another natural choice is $\mat{P} := \mat{B}^{-1}$,
where $\mat{B}$ is a \textit{curvature matrix} of our problem\footnote{See the definition of $\mat{B}$ in our Assumption~\ref{BDef} and
	the corresponding Examples~\ref{ExampleQuadratic}, \ref{ExampleSubspace}, \ref{ExampleSeparable} of different problems.}, which is directly available for the algorithm. 
That resembles the Newton-type direction, 
	and the method with this preconditioner tends to converge
	much faster in practice (see Figure~\ref{fig:Binv}).
However, computing $\mat{B}^{-1}$ (or solving the corresponding linear system with $\mat{B}$) is a very expensive operation in the large scale setting.

Instead of using $\mat{B}^{-1}$, we propose a new family of \textit{Symmetric Polynomial Preconditioners},
that provably improve the spectrum of the objective. The first member of our family is
\beq \label{IntroP1}
\ba{rcl}
\mat{P} & := & \tr(\mat{B}) \mat{I} - \mat{B}.
\ea
\eeq

\begin{figure}[h]
	\centering
	\includegraphics[width=0.70\linewidth]{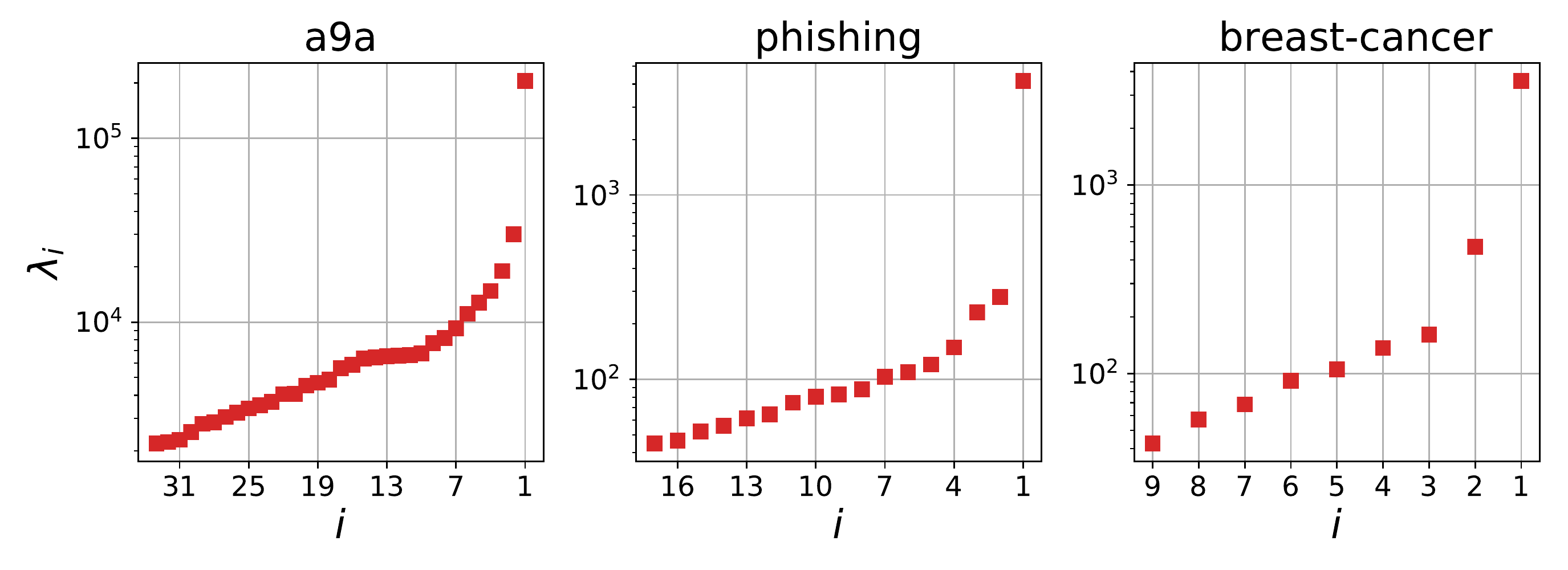}\\[-5pt]
	\includegraphics[width=0.70\linewidth]{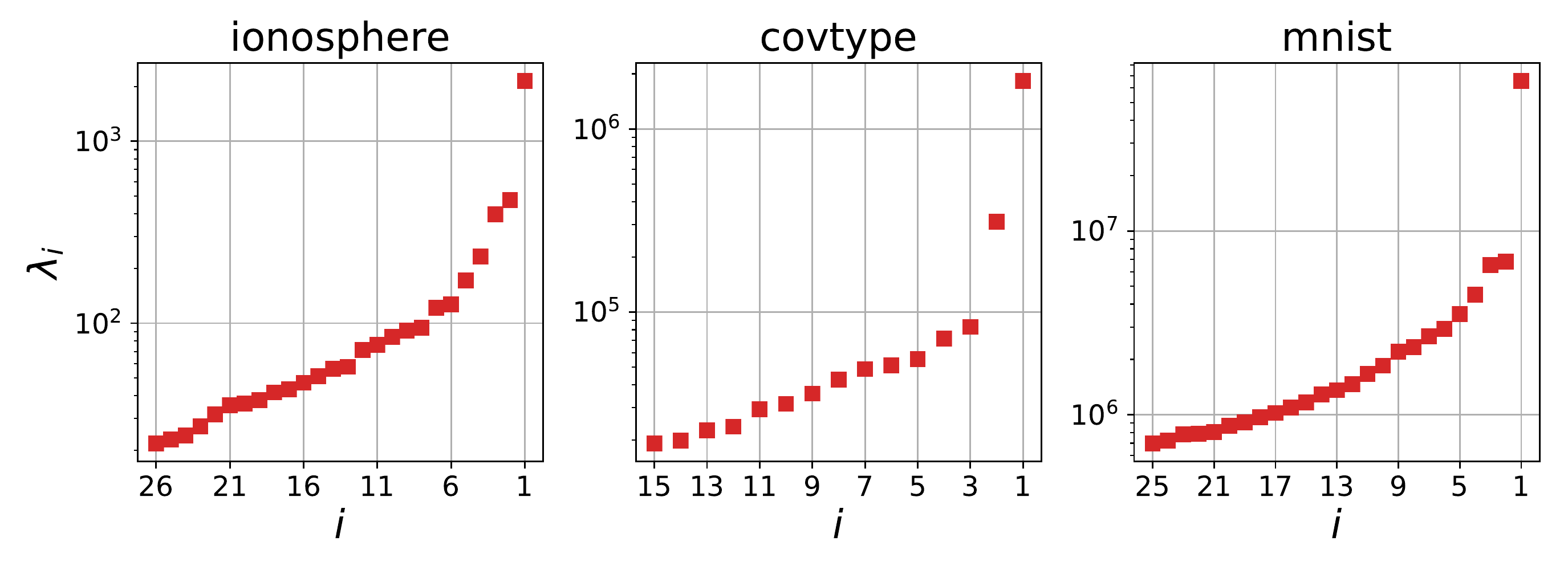}

	\caption{\small Leading eigenvalues (in the logarithmic scale) of the curvature matrix $\mat{B}$,
		for several typical datasets\protect\footnotemark. There are large gaps between the top eigenvalues.}
	\label{fig:Specs}
\end{figure}

We prove that using preconditioner~\eqref{IntroP1} within method \eqref{IntroGM},
makes the condition number \textit{insensitive to the gap between the top two eigenvalues}
of the curvature matrix. Since it is quite common for real data to have a highly nonuniform spectrum
with several large gaps between the top eigenvalues (see Figure~\ref{fig:Specs}),
our preconditioning can significantly improve the convergence of the first-order methods.
At the same time, one step of the form \eqref{IntroGM},\eqref{IntroP1}
is still cheap to compute.
It involves just the standard matrix operations
(trace and the matrix-vector product),
without the need to solve linear systems with the curvature matrix as in Newton's Method.

\footnotetext{\url{https://www.csie.ntu.edu.tw/~cjlin/libsvmtools/datasets/}}

This approach works for general structured nonlinear problems (\textit{not necessarily quadratics})
and also for the problems with possible \textit{composite} parts (e.g., constrained minimization or non-smooth regularization).

Our new \textit{family} of Symmetric Polynomial Preconditioners gradually interpolate between
the first preconditioner~\eqref{IntroP1} and $\mat{P} \propto \mat{B}^{-1}$ as the other extreme case.
We show that increasing the order of the preconditioner, we are able to cut off
\textit{several} top eigenvalues of the curvature matrix, without knowing the actual spectrum.
We can incorporate these preconditioners both into the Gradient Method, as well
as into the accelerated Fast Gradient Method
\cite{nesterov1983method}, with a further provable improvement of the condition number.

Finally, we address the common question of choosing the best possible
preconditioner. We propose a new adaptive strategy for the basic nonlinear Gradient Method
based on the \textit{Krylov subspace} minimization. In this approach,
preconditioner $\mat{P}$ is defined as a general polynomial of the curvature matrix $\mat{B}$
of a fixed (small) degree $\tau$:
$$
\ba{rcl}
\mat{P} & := & a_0 \mat{I} + a_1 \mat{B} + \ldots + a_{\tau} \mat{B}^{\tau},
\ea
$$
where the vector of coefficients $\aa \in \R^{\tau + 1}$ is found by solving
a certain linear system of size $\tau + 1$ in each iteration of the method.
It has a plain interpretation of projecting
the direction $\mat{B}^{-1} \nabla f(\xx_k)$ onto an affine set $\mathcal{K}_{\xx_k}^{\tau}$,
which is the Krylov subspace:
\beq \label{Krylov}
\ba{rcl}
\!\!\!\!
\mathcal{K}_{\xx}^{\tau} & \!\! \Def \!\! &
\LinearHull\bigl\{ \nabla f(\xx), \, \mat{B} \nabla f(\xx), \ldots, \mat{B}^{\tau} \nabla f(\xx)  \bigr\}.
\ea
\eeq
In case of small $\tau$, we can solve this linear system easily
and obtain the best preconditioning guarantee for our method,
which is \textit{adaptive} for each iteration.

\paragraph{Related Work.}

It is widely known that
the standard Conjugate Gradient Method
is \textit{optimal}
in the class of the first-order algorithms
for unconstrained minimization of convex \emph{quadratic} functions
\cite{nemirovski1995information}.
The $k$th iteration of the Conjugate Gradients finds
the full minimum of the objective over the $k$-dimensional Krylov subspace,
and thus it provably solves the problem after $k = n$ iterations,
where $n$ is the dimension of the problem.
Quadratic minimization is equivalent to solving a system of linear equations,
therefore it is often referred as the \textit{linear} case.
Polynomial preconditioning 
for solving large linear systems has been extensively studied 
during the last several decades;
see \cite{dubois1979approximating,johnson1983polynomial,saad1985practical,van1995polynomial,liu2015polynomial,loe2022toward}
and references therein.
See also Section~\ref{SubsectionDiscussion} for the comparison
of our preconditioning strategies with
the linear Conjugate Gradient Method.

The situtation with \textit{nonlinear} problems is more difficult.
Along with the basic Gradient Method,
the classical approaches include
the Nonlinear Conjugate Gradients
and Quasi-Newton Methods (see, e.g. \cite{nocedal2006numerical}),
which typically demonstrate a decent practical performance, 
while replicating the standard Conjugate Gradients
in the linear case. However, these methods \textit{lack} of having
any good \textit{global complexity bounds}, and thus 
in the worst-case scenario they can actually perform even worse than the Gradient Method
\cite{gupta2023nonlinear}.
At the same time, the Fast Gradient Method developed by \cite{nesterov1983method}
is \textit{optimal} for the class of nonlinear problems
with a \textit{uniformly bounded eigenvalues} of the Hessian \cite{nemirovskii1983problem}.
This assumption does not take into account the actual distribution of the spectrum.
Hence, it can not distinguish the problems with large gaps between the top eigenvalues, 
as in Figure~\ref{fig:Specs}.

There have been several attempts to study more specific problem formulations,
and so to gain a provable advantage
for the optimization algorithms
by leveraging the spectrum of the Hessian.
Thus,
the quadratic minimization problems were studied under
the assumption of a particular 
\textit{probability distribution} for the eigenvalues
\cite{scieur2020universal,cunha2022only},
or assuming a certain fixed \textit{spectral gap} \cite{goujaud2022super},
revealing the advantages of employing the Heavy-ball Method \cite{polyak1987introduction}
in these cases.
Another example is the Stochastic Spectral Descent \cite{kovalev2018stochastic},
which improves the condition number 
for quadratic problems if we know some
of the eigenvectors.

In this work, we consider a refined smoothness characterization 
of the objective
with the curvature matrix $\mat{B}$
(Assumption~\ref{BDef}).
It is similar in spirit to that one used
in 
Stochastic Dual Newton Ascent
\cite{qu2016sdna}.
An important particular instance of this class of algorithms
is the Randomized Coordinate Descent with \textit{Volume Sampling}
\cite{rodomanov2020randomized}. In the latter method,
it was proposed to select subsets of variables of
certain size $m$ proportionally to the determinants
of principal submatrices of $\mat{B}$.
While this approach was practically implementable
only for the subsets of size $m = 1$ or $2$,
it was shown that, in theory, the method is insensitive 
to the large spectral gap between the top $m - 1$
eigenvalues.

Surprisingly, our new family of the Symmetric Polynomial Preconditioners
can be viewed as a \textit{deterministic version}
of the Volume Sampling technique (with $m = \tau + 1$
where $\tau$ is the degree of a preconditioning polynomial;
preconditioner \eqref{IntroP1} corresponds to $\tau = 1$).
Thus, we provide the Volume Sampling with a novel deterministic interpretation, which also leads
to new accelerated and composite optimization algorithms
(see Section~\ref{SubsectionStochastic} for a detailed comparison).


\paragraph{Contributions.}

\begin{table*}
	\centering
	\small
	\renewcommand{\arraystretch}{1.4}

	\newcolumntype{\CenteredParagraphColumn}[1]{>{\centering\arraybackslash}m{#1}}
	\newcolumntype{\Strategy}{>{\raggedright\arraybackslash}m{0.28\linewidth}}
	\newcolumntype{\Preconditioner}{m{0.15\linewidth}}
	\newcolumntype{\ConditionNumber}{\CenteredParagraphColumn{0.2\linewidth}}
	\newcolumntype{\Methods}{\CenteredParagraphColumn{0.08\linewidth}}
	\newcolumntype{\Cost}{\CenteredParagraphColumn{0.13\linewidth}}

	\begin{tabular}{\Strategy @{\hspace{1.5em}} \Preconditioner \ConditionNumber \Methods \Cost}
		 & Preconditioner & Cond.\ number, $\beta / \alpha$ & Methods & Cost \\
		\toprule
		Classic Gradient Method &  $\mat{P} = \mat{I}$  &  $\lambda_1 / \lambda_n$
		& GM, FGM & \GoodProperty{cheap} \\
		\hline
		``Full Preconditioning'' &  $\mat{P} = \mat{B}^{-1}$ & $1$
		& GM, FGM & \BadProperty{expensive} \\
		\midrule
		Symmetric Polynomial\linebreak Preconditioning \GoodProperty{(ours)} &
		$\mat{P} = \mat{P}_{\tau}$ \eqref{PRecDef} & $\lambda_1 / \lambda_n \cdot \xi_{\tau}(\blambda)$
		& GM, FGM & \GoodProperty{cheap\linebreak for small $\tau$} \\
		\hline
		Krylov Subspace\linebreak Minimization \GoodProperty{(ours)} &
		optimal poly. & $\lambda_{\tau + 1} / \lambda_n$ &
		GM & \GoodProperty{cheap\linebreak for small $\tau$}  \\
		\bottomrule
	\end{tabular}
	\caption{$\frac{\beta}{\alpha}$
		for different preconditioning strategies, $\blambda = \blambda(\mat{B})$.
		Note that $\xi_{\tau}(\blambda) \leq 1$, and $\xi_{\tau}(\blambda) \to 0$
		in case of large spectral gaps, namely
		when $\frac{\lambda_{1}}{\lambda_{\tau + 1}} \to \infty$
		(see Section~\ref{SectionSymmetric}).
		For solving the problem
		with $\epsilon$-accuracy,
		GM needs $k(\epsilon) = \cO( \frac{\beta}{\alpha} \cdot \frac{1}{\epsilon} )$
		and
		$k(\epsilon) = \cO( \frac{\beta}{\alpha} \cdot \frac{L}{\mu} \cdot \log \frac{1}{\epsilon}  )$
		iterations for convex and strongly convex functions, respectively. FGM needs only $\sqrt{k(\epsilon)}$ iterations
	    (Theorems~\ref{TheoremGM} and~\ref{TheoremFGM}).}
	\label{TableComplexities}
\end{table*}

We propose several polynomial preconditioning strategies 
for first-order methods for solving
a general composite convex optimization problem,
and prove their better global complexity guarantees,
specifically:

\begin{itemize}
	
	\item We study the convergence of the basic Gradient Method (GM, Algorithm~\ref{alg:GM})
	and the accelerated Fast Gradient Method (FGM, Algorithm~\ref{alg:FGM}) with a general
	(arbitrarily fixed) preconditioning matrix. We introduce \textit{two} condition numbers,
	that are designated to the different parts of the objective ($L / \mu$ for nonlinearity and 
	$\beta / \alpha$ for the curvature matrix),
	and show that they serve as main complexity factors.

	\item We develop a new family of Symmetric Polynomial Preconditioners (Section \ref{SectionSymmetric}).
	Combining them with the preconditioned Gradient Methods,
	we establish a significant improvement of the curvature condition number $\beta / \alpha$
	in case of \textit{large gaps} between the top eigenvalues
	of the matrix (see Table~\ref{TableComplexities}).
	
	\item Then,  
	we propose a new adaptive procedure based on the Krylov subspace minimization
	(Algorithm~\ref{alg:KrylovGM})
	that achieves the \textit{best polynomial} preconditioning.
	We present the guarantees we can get,
	including cutting off the top eigenvalues directly and by employing the Chebyshev polynomials,
	and compare this approach with the Symmetric Polynomial Preconditioning.
	
	\item Numerical experiments are provided. 
	
\end{itemize}


  \section{Notation and Assumptions}
\label{SectionNotation}

We consider the following optimization problem
given in the \textit{composite} form:
\beq \label{MainProblem}
\ba{rcl}
F^{\star}
\;\; = \;\;
\min\limits_{\xx \in \R^n} \Bigl\{
F(\xx) & \Def & f(\xx) + \psi(\xx)
\Bigr\},
\ea
\eeq
where $f: \R^n \to \R$ is a differentiable convex function which
is the \textit{main} part of the problem,
and $\psi: \R^n \to \R \cup \{ +\infty \}$
is a proper closed convex function that can be nondifferentiable
but has a \textit{simple} structure.
For example, it can be an indicator of a
convex set, or $\ell_1$-regularizer.

Additionally, we fix some symmetric positive-definite matrix
$\mat{B} \in \R^{n \times n}$ (notation $\mat{B} = \mat{B}^\top \succ 0$).
This matrix plays the key role in our characterization of the smoothness
properties of $f$. Namely, we assume the following
(considering for simplicity two-times differentiable functions):

\begin{framedassumption} \label{BDef}
  The Hessian of $f$ is uniformly bounded,
  for some constants $L \geq \mu \geq 0$:
  \beq \label{HessBound}
  \ba{rcl}
  \mu \mat{B} & \preceq & \nabla^2 f(\xx)
  \;\; \preceq \;\;
  L \mat{B}, \qquad \forall \xx \in \R^n.
  \ea
  \eeq
\end{framedassumption}

Having fixed the matrix $\mat{B}$,
we define the corresponding induced norm by
$ \|\xx\|_{\mat{B}}  \Def  \la \mat{B} \xx, \xx \ra^{1/2}$,
$\xx \in \R^n$.
Thus, matrix $\mat{B}$ is responsible for fixing the coordinate system in the problem.
Then, condition \eqref{HessBound} can be rewritten
in terms of the global lower and upper bound on the
first-order approximation of $f$ \cite{nesterov2018lectures}:
\vspace{-1mm}
\beq \label{FuncGlobalBounds}
\ba{rcl}
\!\!\!
\frac{\mu}{2}\|\yy - \xx\|_{\mat{B}}^2
& \leq &
f(\yy) - f(\xx) - \la \nabla f(\xx), \yy - \xx \ra \\[10pt]
& \leq & \frac{L}{2}\|\yy - \xx\|_{\mat{B}}^2,
\qquad \forall \xx, \yy \in \R^n.
\ea
\eeq

  In what follows, we denote by $\blambda = \blambda(\mat{B}) \in \R^n$
the vector of eigenvalues for the matrix $\mat{B}$,
sorted in a nonincreasing order:
$
\lambda_1  \, \geq \, \lambda_2
\, \geq \, \ldots \, \geq \, \lambda_n$.


The classical example is $\mat{B} := \mat{I}$ (identity matrix).
Then, condition \eqref{HessBound} implies that
the function $f$ is (strongly) convex and has the Lipschitz continuous
gradient.
However, by choosing a specific $\mat{B}$,
we tend to achieve a better granularity of the description of our problem class
and thus to improve the convergence properties of the methods.

\begin{example} \label{ExampleQuadratic}
  Let $\aa \in \R^{n}$. Then, the quadratic function
  \vspace{-1mm}
  $$
  \ba{rcl}
  f(\xx) & = & \frac{1}{2}\la \mat{B} \xx, \xx \ra
  - \la \aa, \xx \ra,
  \ea
  $$
  satisfies condition \eqref{HessBound} with
  $L = \mu = 1$.
\end{example}

We see that in this case, the so-called \textit{condition number}
$L / \mu$ is just $1$, which means that preconditioning
the Gradient Method \eqref{IntroGM} with
the matrix $\mat{P} :=\mat{B}^{-1}$ would give an immediate
convergence to the solution.
However, inverting the matrix is prohibitively expensive for large scale problems.
Our aim is to find a suitable trade-off
between improving the condition number
and the arithmetic cost of algorithm steps.
Let us consider the following important example
which can be met in many practical applications.

\begin{example} \label{ExampleSubspace}
  Let $\mat{A} \in \R^{m \times n}$ be a given data matrix,
  and $\bb \in \R^m$ be a given vector. Denote,
  \vspace{-1mm}
  $$
  \ba{rcl}
  f(\xx) &=& g(\mat{A} \xx + \bb)
  \ea
  $$
  Then, the derivatives are as follows:
  $\nabla f(\xx) = \mat{A}^\top \nabla g(\mat{A} \xx + \bb)$ and
  $\nabla^2 f(\xx)  = \mat{A}^\top \nabla^2 g(\mat{A} \xx + \bb) \mat{A}$.
  Hence, assuming:
  $ \mu \mat{I_m} \, \preceq \, \nabla^2 g(\xx)
  \, \preceq \, L \mat{I_m}$, $\forall \xx$, with some $L \geq \mu \geq 0$,
  condition \eqref{HessBound} is satisfied~\footnote{Here, we assume that
    $\mat{A}^\top \mat{A} \succ 0$ which is typically the case when $m \ggg n$. Otherwise, we can
    reduce the dimensionality.} with
  \vspace{-1mm}
  $$
  \ba{rcl}
  \mat{B} & := & \mat{A}^\top \mat{A}.
  \ea
  $$
  At the same time, for $\mat{B} := \mat{I_n}$
  (the standard Euclidean norm), the Lipschitz constant
  increases by the factor $\lambda_1(\mat{A}^{\top} \mat{A})$,
  which makes the problem extremely ill-conditioned.
\end{example}
A particular case of this example is \textit{separable optimization},
or \textit{generalized linear models} \cite{bishop2006pattern},
which covers the classical regression and classification models.
\begin{example} \label{ExampleSeparable}
  Let
  \vspace{-1mm}
  $$
  \ba{rcl}
  f(\xx) & = & \frac{1}{m} \sum\limits_{i = 1}^m
  \phi(\la \aa_i, \xx \ra),
  \qquad \xx \in \R^n,
  \ea
  $$
  where $\phi : \R \to \R$ is a loss function satisfying:
  $\mu \leq \phi''(t) \leq L$,  $\forall t \in \R$,
  with some $L \geq \mu \geq 0$.
  Then, forming the matrix $\mat{A} \in \R^{m \times n}$
  whose rows are $\aa_1^\top, \ldots, \aa_m^\top$
  and setting $\mat{B} := \mat{A}^\top \mat{A},$
  condition \eqref{HessBound} holds.
\end{example}

  \section{Preconditioned Gradient Methods}
\label{SectionGM}

A natural intention would be to use the global upper bound
\eqref{FuncGlobalBounds} as a model for the smooth part of the objective.
However, the direct minimization of such upper model
requires to solve the linear system with the matrix $\mat{B}$,
which can computationally unfeasible
for large scale problems.

Instead, let us fix for our \textit{preconditioner} some positive definite
symmetric matrix $\mat{P} = \mat{P}^{\top} \succ 0$,
which satisfies the following bound,
for some $\alpha := \alpha(\mat{P})$ and $\beta := \beta(\mat{P}) \geq \alpha > 0$:
\beq \label{PBound}
\boxed{
  \ba{rcl}
  \alpha \mat{B}^{-1} & \preceq &
  \mat{P} \;\; \preceq \;\; \beta \mat{B}^{-1}.
  \ea
}
\eeq
We are going to use this matrix instead of~$\mat{B}^{-1}$ in our methods.
For a fixed symmetric positive definite matrix $\mat{P}$ and parameter $M > 0$,
we denote the \textit{gradient step} from a point $\xx \in \dom \psi$ along a
\textit{gradient direction} $\Vector{g} \in \R^n$
by
\begin{align*}
  \hspace{2em}&\hspace{-2em}
  \GradientStep_{M, \mat{P}}(\xx, \Vector{g})
  \\
  &\Def
  \argmin_{\yy \in \dom \psi}\Bigl\{
  \la \Vector{g}, \yy \ra + \psi(\yy)
  + \frac{M}{2}\|\yy - \xx\|_{\mat{P}^{-1}}^2 \Bigr\}.
\end{align*}
This operation is well-defined since the objective function in the above
minimization problem is strongly convex.
We assume that both $\psi$ and $\mat{P}$ are reasonably simple so that
the corresponding gradient step can be efficiently computed.
An important case is $\psi = 0$ for which we have
$
  \GradientStep_{M, \mat{P}}(\xx, \Vector{g})
  =
  \xx - \frac{1}{M} \mat{P} \Vector{g}
$.
The latter expression can be efficiently computed whenever one can cheaply
multiply the matrix~$\mat{P}$ by any vector.

\subsection{Preconditioned Basic Gradient Method}
\label{SectionBasicGM}

First, we consider the basic first-order scheme shown in \cref{alg:GM}
for solving the composite problem~\eqref{MainProblem}.
For simplicity, in this section, we only present a version of this method
with a fixed step size and assume that all necessary constants are known.
An adaptive version of \cref{alg:GM} which does not have these limitations
and is more efficient in practice can be found in \cref{SectionAdaptive}.

\newpage

\begin{algorithm}[h!]
  \caption{Preconditioned Basic Gradient Method}
  \label{alg:GM}
  \begin{algorithmic}
    \STATE {\bfseries Input:} $\Vector{x}_0 \in \dom \psi$,
    $\mat{P} = \mat{P}^{\top} \succ 0$, $M > 0$.
    \FOR{$k = 0, 1, \dots$}
    \STATE Compute
    $
      \Vector{x}_{k + 1}
      =
      \GradientStep_{M, \mat{P}}\bigl( \Vector{x}_k, \nabla f(\Vector{x}_k) \bigr)
    $.
    \ENDFOR
  \end{algorithmic}
\end{algorithm}

For \cref{alg:GM}, we can prove the following results.

\begin{framedtheorem} \label{TheoremGM}
  Consider \cref{alg:GM} with $M = \beta L$.
  Then, at each iteration $k \geq 1$, we have
  \begin{equation}
    \label{Alg1Convex}
    F(\Vector{x}_k) - F^{\star}
     \;\; \leq \;\;
    \frac{\beta}{\alpha} 
    \frac{L \RelativeNorm{\Vector{x}_0 - \Vector{x}^{\star}}{\mat{B}}^2}{k}.
  \end{equation}
  When $\mu > 0$, the convergence is linear: for all $k \geq 1$,
  \begin{equation}
    \label{Alg1Strongly}
    F(\Vector{x}_k) - F^{\star}
    \;\; \leq \;\;
    \Bigl( 1 - \frac{1}{4} \frac{\alpha}{\beta} \frac{\mu}{L} \Bigr)^k
    [F(\Vector{x}_0) - F^{\star}].
  \end{equation}
\end{framedtheorem}

We see that one of the principal complexity factors in the above estimates
is the condition number $\beta / \alpha$ which depends on the choice
of our preconditioner $\mat{P}$ (see~\eqref{PBound}).
For the basic choice $\mat{P} = \mat{I}$, we have
$\beta / \alpha = \lambda_1 / \lambda_n$.
However, as we show in the following sections, it is possible to use
more efficient (and still quite cheap) preconditioners which improve this
condition number.

\subsection{Preconditioned Fast Gradient Method}
\label{SectionFastGM}

Now let us consider an accelerated scheme shown in \cref{alg:FGM}.
This algorithm is one of the standard variants of the Fast Gradient Method (FGM)
known as the Method of Similar Triangles
(see, e.g., Section~6.1.3 in~\cite{nesterov2018lectures})
but adapted to our assumptions~\eqref{HessBound} and~\eqref{PBound}.

\begin{algorithm}
  \caption{Preconditioned Fast Gradient Method}
  \label{alg:FGM}
  \begin{algorithmic}
    \STATE {\bfseries Input:} $\Vector{x}_0 \in \dom \psi$,
    $\mat{P} = \mat{P}^{\top} \succ 0$, $M > 0$, $\rho \geq 0$.
    \STATE Set $\Vector{v}_0 = \Vector{x}_0$, $A_0 = 0$.
    \FOR{$k = 0, 1, \dots$}
    \STATE
      Find $a_{k + 1}$ from eq.\
      $\frac{M a_{k + 1}^2}{A_k + a_{k + 1}} = 1 + \rho (A_k + a_{k + 1})$.
    \STATE
      Set $A_{k + 1} = A_k + a_{k + 1}$,
      $H_k = \frac{1 + \rho A_{k + 1}}{a_{k + 1}}$.
    \STATE
      Set $\theta_k = \frac{a_{k + 1}}{A_{k + 1}}$,
      $\omega_k = \frac{\rho}{H_k}$,
      $\gamma_k = \frac{\omega_k (1 - \theta_k)}{1 - \omega_k \theta_k}$.
    \STATE
      Set
      $\hat{\Vector{v}}_k = (1 - \gamma_k) \Vector{v}_k + \gamma_k \Vector{x}_k$.
    \STATE
      Set $\Vector{y}_k = (1 - \theta_k) \Vector{x}_k + \theta_k \hat{\Vector{v}}_k$.
    \STATE
      Compute
      $
        \Vector{v}_{k + 1}
        =
        \GradientStep_{H_k, \mat{P}}\bigl(
          \hat{\Vector{v}}_k, \nabla f(\Vector{y}_k)
        \bigr)
      $.
    \STATE
      Set
      $
        \Vector{x}_{k + 1}
        =
        (1 - \theta_k) \Vector{x}_k + \theta_k \Vector{v}_{k + 1}
      $.
    \ENDFOR
  \end{algorithmic}
\end{algorithm}

As in other versions of FGM, to properly handle strongly convex problems,
\cref{alg:FGM} requires the knowledge of the strong convexity parameter $\alpha$ and $\mu$
(or, more precisely, their product $\rho = \alpha \mu$).
For non-strongly convex problems, we can always choose $\alpha = \mu = 0$.
See also \cref{SectionAdaptive} for a variant of \cref{alg:FGM} which
can automatically adjust the constant~$M$ in iterations.

The convergence results for \cref{alg:FGM} are as follows.

\begin{framedtheorem} \label{TheoremFGM}
  Consider \cref{alg:FGM} with $M = \beta L$
  and $\rho = \alpha \mu$.
  Then, at each iteration $k \geq 1$, we have
  \begin{equation}
    \label{Alg2Convex}
    F(\Vector{x}_k) - F^{\star}
    \;\; \leq \;\;
    2 \,
   \frac{\beta}{\alpha} 
    \frac{L \RelativeNorm{\Vector{x}_0 - \Vector{x}^{\star}}{\mat{B}}^2}{k^2}.
  \end{equation}
  When $\mu > 0$, the convergence is linear: for all $k \geq 1$,
  \[
    F(\Vector{x}_k) - F^{\star}
    \;\; \leq \;\;
    \biggl( 1 - \sqrt{\frac{\alpha}{\beta} \frac{\mu}{L}} \, \biggr)^{k - 1}
    \frac{\beta}{\alpha}
    \frac{L}{2} \RelativeNorm{\Vector{x}_0 - \Vector{x}^{\star}}{\mat{B}}^2.
  \]
\end{framedtheorem}

Comparing these estimates with those from \cref{TheoremGM},
we see that the accelerated scheme is much more efficient.
For instance, to reach accuracy $\epsilon > 0$ in terms of the objective
function in the non-strongly convex case, \cref{alg:GM} needs
$
  k(\epsilon)
  =
  \frac{\beta}{\alpha}
  \frac{L \RelativeNorm{\Vector{x_0} - \Vector{x}^{\star}}{\mat{B}}^2}{\epsilon}
$
iterations, while for \cref{alg:FGM} this number is only
$k_2(\epsilon) = \sqrt{2 k(\epsilon)}$.
Similar conclusions are valid in the strongly convex case.

Despite having much weaker dependency on the condition number $\beta / \alpha$,
\cref{alg:FGM} is still quite sensitive to it.
Thus, the proper choice of the preconditioner~$\mat{P}$ is important for both
our methods.

  \section{Symmetric Polynomial Preconditioning}
\label{SectionSymmetric}

We would like to have a \textit{family} of preconditioners
$\mat{P}_{\tau}$ for our problem indexed by some parameter $\tau$.
Varying $\tau$ should provide us with  a trade off between
the spectral quality of approximation \eqref{PBound} of the inverse matrix
and the arithmetical cost of computing the preconditioner.

Surprisingly, such a family of preconditioners can be built by using
\textit{symmetric polynomials}, the classical objects of Algebra.
We prove that our preconditioning improves the condition number $\frac{\beta}{\alpha}$
of the problem, by automatically
cutting off the large gaps between the top eigenvalues.

\subsection{Definition and Basic Properties}
\label{sec:DefinitionOfPreconditioners}

We define
the family of symmetric matrices $\{ \mat{P}_{\tau} \}_{0 \leq \tau \leq n - 1}$
recursively. We start with identity matrix:
$ \mat{P}_0 \; \Def \; \mat{I} $.
Then, 
\vspace{-1mm}
\beq \label{PRecDef}
\boxed{
\ba{rcl}
\mat{P}_{\tau} & \Def &
\frac{1}{\tau}
\sum\limits_{i = 1}^{\tau} (-1)^{i - 1}
\mat{P}_{\tau - i} \mat{U}_{i},
\ea
}
\eeq
where $\mat{U}_\tau  \Def \tr(\mat{B}^{\tau}) \mat{I} - \mat{B}^{\tau}$
are the auxiliary matrices.
It turns out that matrices \eqref{PRecDef}
serve as a good approximation of the inverse matrix:
$ \mat{P}_{\tau} \approx  \mat{B}^{-1}$,
up to some multiplicative constant, and the quality of such approximation
is gradually improving when increasing parameter $\tau$.
Let us look at several first members. Clearly,
\vspace{-1mm}
\beq \label{P1Def}
\ba{rcl}
\mat{P}_{1} & = & \tr(\mat{B}) \mat{I} - \mat{B},
\ea
\eeq
which is very easy to handle, by having computed the trace of the curvature matrix. Then, multiplying $\mat{P}_1$
by any vector would require just one matrix-vector multiplication with
our original $\mat{B}$.
Further,
\vspace{-1mm}
\begin{equation} \label{P2Def}
  \begin{array}{rcl}
    \mat{P}_{2} &=& \tfrac{1}{2} \tr(\mat{P}_1 \mat{B})
    \mat{I} - \mat{P}_1 \mat{B}
    \\
    \\
    &=&
    \tfrac{1}{2} \bigl[ [\tr(\mat{B})]^2 - \tr(\mat{B}^2) \bigr] \mat{I} - \tr({\mat{B}}) \mat{B} + \mat{B}^2,
  \end{array}
\end{equation}
thus its use would cost just \textit{two} matrix-vector products with $\mat{B}$,
having evaluated\footnote{Note
  that $\tr(\mat{B}^2) = \sum_{i = 1}^n \| \mat{B}[:, i] \|_2^2$, where $\mat{B}[:, i] \in \R^n$
  is the $i$th column of $\mat{B}$.}${}^{,}$\footnote{%
  For general $\tau$, we can also use a stochastic estimate of the trace: $\xi_{\tau} \Def n \la \mat{B}^{\tau} \uu, \uu \ra$,
  where $\uu \in \R^n$ is uniformly distributed on the unit sphere. It would give an unbiased estimate:
  $\mathbb{E}[\xi_{\tau}] = n \mathbb{E}[ \tr (\uu^{\top} \mat{B}^\tau \uu) ] =
  n \tr (\mathbb{E}[ \uu \uu^{\top} ] \mat{B}^{\tau}) = \tr(\mat{B}^{\tau})$.
} the numbers $\tr(\mat{B})$ and $\tr(\mat{B}^2)$.

It is clear that in general $\mat{P}_{\tau} = p_\tau(\mat{B})$, where $p_\tau$ is a polynomial of a fixed degree $\tau$
with coefficients that can be found recursively from \eqref{PRecDef}.
Let us give a useful interpretation for our family of preconditioners, that also explains their name.
For $\aa \in \R^{n - 1}$, we denote by $\sigma_0(\aa), \ldots, \sigma_{n - 1}(\aa)$
the \textit{elementary symmetric polynomials}
in $n - 1$ variables\footnote{That is
  $
  \sigma_{\tau}(\aa) \; \Def \; \sum_{1 \leq i_1 < \ldots < i_{\tau} \leq n-1}
  a_{i_1} \cdot \ldots a_{i_{\tau}}.
  $}.
It is known that every symmetric polynomial (that is invariant to any permutation of the variables) can be represented as a weighed sum of elementary symmetric polynomials \cite{dummit2004abstract}.
We establish the following important characterization.

\begin{lemma} \label{LemmaSpec}
  Let $\mat{B} = \mat{Q} \Diag( \blambda ) \mat{Q}^{\top}$ be the spectral decomposition. Then,
  \beq \label{P_tau_spec}
  \ba{rcl}
  \!\!\!\! \mat{P}_{\tau}
  & \!\!\! = \!\!\! &
  \mat{Q} \Diag( \sigma_{\tau}(\blambda_{-1}),
  \ldots,
  \sigma_{\tau}(\blambda_{-n})   ) \mat{Q}^{\top},
  \ea
  \eeq
  where $\blambda_{-i} \in \R^{n - 1}$ is the vector that
  contains all elements of $\blambda$ except $\lambda_i$.
\end{lemma}

In particular, we justify $\mat{P}_{\tau} \succ 0$.
For $\tau = n - 1$, we get
\vspace{-2mm}
\beq \label{PN_1}
\ba{rcl}
\mat{P}_{n - 1} & \overset{\eqref{LemmaSpec}}{=} &  \det(\mat{B}) \mat{B}^{-1}
\;\,\Def \;\, \Adj(\mat{B}),
\ea
\eeq
which gives us the true inverse matrix $\mat{B}^{-1}$
up to the constant factor $\det(\mat{B})$.

\newpage
\subsection{Approximation Quality}
\label{SubsectionApproxQuality}

Now, let us show that the quality of approximation $\mat{P}_{\tau} \approx \mat{B}^{-1}$
and that the corresponding condition number $\frac{\beta}{\alpha}$
is improving when $\tau$ is increasing.

\begin{framedtheorem} \label{TheoremPApprox}
  For any $\tau$, we have
  \beq \label{PApproxQuality}
  \ba{rcl}
  \!\!\!\!\!\!
  \lambda_n \sigma_{\tau}(\blambda_{-n}) \mat{B}^{-1}
  & \!\! \preceq \!\! &
  \mat{P}_{\tau}
  \; \preceq \;
  \lambda_1 \sigma_\tau(\blambda_{-1}) \mat{B}^{-1}
  \ea
  \eeq
  Therefore, the condition number is bounded as
  $$
  \ba{rcl}
  \frac{\beta}{\alpha}
  & \overset{\eqref{PApproxQuality}}{=} &
  \frac{\lambda_1}{\lambda_n}
  \cdot \xi_\tau(\blambda),
  \quad
  \text{where}
  \quad
  \xi_{\tau}(\blambda) \; \Def \;
  \frac{\sigma_{\tau}(\blambda_{-1})}{ \sigma_{\tau}( \blambda_{-n} ) }.
  \ea
  $$
\end{framedtheorem}

\begin{figure}[h!]
	\centering
	\includegraphics[width=0.7\linewidth]{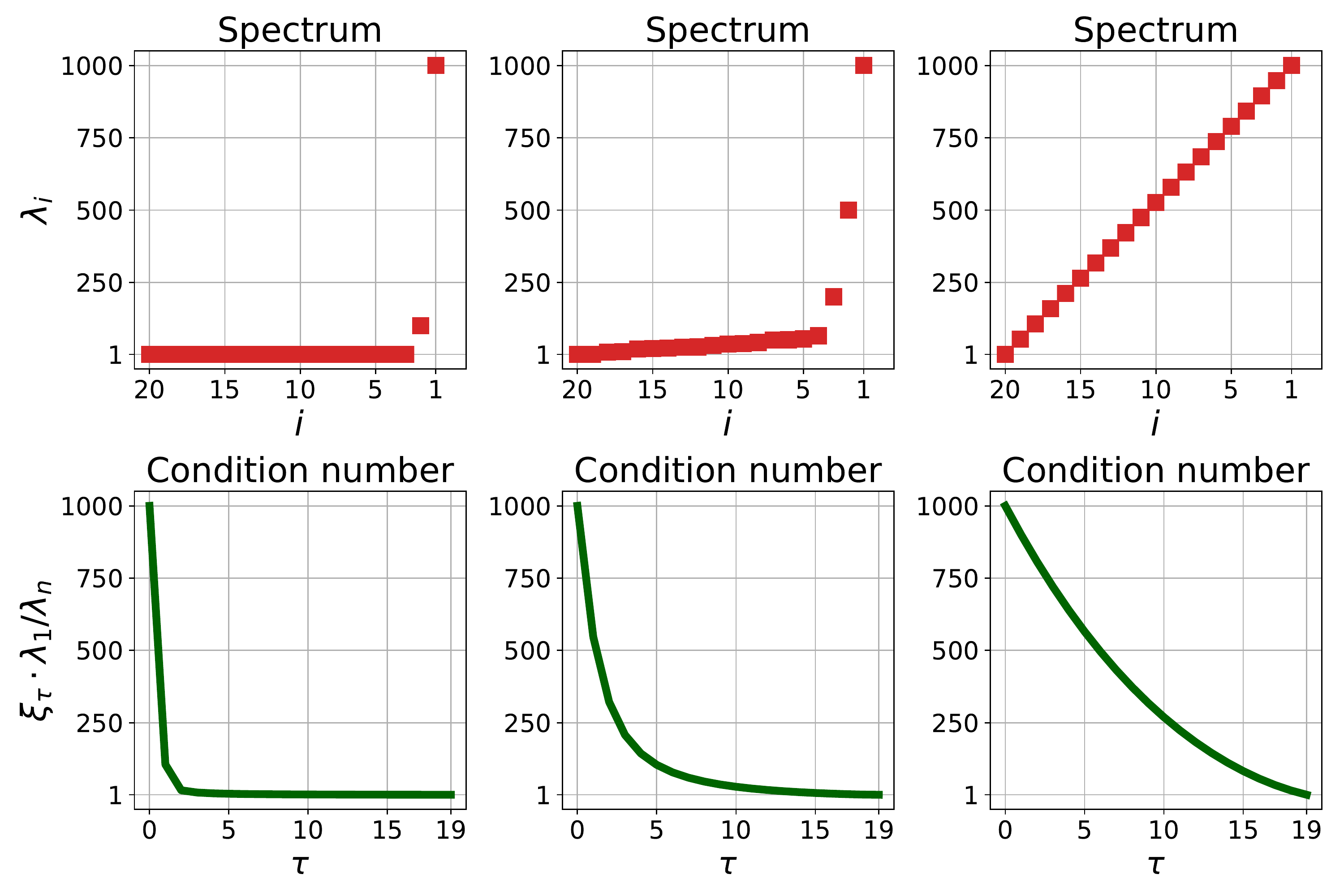}
	\caption{
		\textit{Above:} different distributions of eigenvalues $\blambda(\mat{B})$.
		\textit{Below:} the corresponding improvement of the condition number 
		$\frac{\beta}{\alpha} = \xi_{\tau}(\blambda) \cdot \lambda_1 / \lambda_n$
		by using the preconditioner $\mat{P}_{\tau}$
		of order $\tau$.}
	\label{fig:Xi}
\end{figure}

Note that  $\xi_{\tau}(\blambda) \leq 1$.
It is equal to $1$ for $\tau = 0$ and can be much smaller for bigger values of $\tau$.
For example, for $\tau = 1$ (thus using preconditioner $\mat{P}_1$ given by \eqref{P1Def}), we get
$$
\ba{rcl}
\xi_{1}(\blambda)
& = &
\frac{\sum_{i = 2}^n \lambda_i }{ \sum_{i = 1}^{n - 1} \lambda_i}
\;\; \leq \;\; 1,
\ea
$$
and it is much smaller than $1$ when $\lambda_1 \ggg \lambda_2$,
which corresponds to the case when the highest eigenvalue
is well separated from the others.
Therefore, the methods with preconditioner $\mat{P}_1$ achieve
a \textit{provable acceleration} in the case of large gap between $\lambda_1$
and $\lambda_2$,
\textit{without explicit knowledge} of the spectrum of $\mat{B}$.
The price of using $\mat{P}_1$ instead of $\mat{P}_0$
is just one extra matrix-vector product per iteration.
Let us summarize the main properties of $\xi_{\tau}(\blambda)$
(see also Figure~\ref{fig:Xi}).
\begin{lemma} \label{LemmaXiProperties}
	It holds that $	\xi_{0}(\blambda)  =  1$, $\xi_{n - 1}(\blambda) = \frac{\lambda_n}{\lambda_1}$,
	and $\xi_{\tau}(\blambda)$ monotonically decreases with $\tau$. Moreover,
	$$
	\ba{rcl}
	\xi_{\tau}(\blambda) & \to & 0 \quad \text{when} \quad \frac{\lambda_{1}}{\lambda_{\tau + 1}}
	\; \to \; \infty.
	\ea
	$$
\end{lemma}

We see that  $\tau$ interpolates
$\mat{P}_{\tau}$ between
$\mat{I}$
and $\Adj(\mat{B})$,
while the condition number $\frac{\beta}{\alpha}$ changes from  
$\frac{\lambda_1}{\lambda_n}$ to $1$
Therefore, we obtain an extra degree of freedom in our methods
for choosing an appropriate small value of $\tau$,
that improves the spectrum of the problem
by cutting off large gaps between $\lambda_1$ and $\lambda_{\tau + 1}$.
\subsection{Stochastic Representation}
\label{SubsectionStochastic}

Let us provide another interesting interpretation for our family of
preconditioners.
One way of approximating the inverse matrix $\mat{B}^{-1}$ could be to extract
from~$\mat{B}$ a randomly selected principal submatrix of size~$\tau + 1$,
compute its inverse, put it back into the original ``big matrix''
and zero out all other elements outside the submatrix.
It turns out that, in expectation, the result of this operation is exactly
proportional to our preconditioner~$\mat{P}_{\tau}$ if we pick the submatrix
from a special \emph{volume sampling}
distribution~\cite{Deshpande.etal-06-MatrixApproximation}.

\begin{theorem} \label{TheoremVolumeSampling}
  For any $0 \leq \tau \leq n - 1$, it holds that
  \begin{equation}
    \label{PrecondVolumeSampling}
    \mat{P}_{\tau}
    \propto
    \Expectation_{S \DistributedAs \VolumeSampling_{\tau + 1}(\mat{B})} [
      \mat{I}_S (\mat{B}_{S \times S})^{-1} \mat{I}_S\Transpose
    ],
  \end{equation}
  where $S \subseteq \Set{1, \ldots, n}$ is a random $(\tau + 1)$-element
  subset of coordinates,
  $\mat{I}_S \in \RealMatrices{n}{(\tau + 1)}$ is the matrix obtained from the
  identity matrix by retaining only the columns with indices from~$S$,
  $
    \mat{B}_{S \times S}
    \Def
    \mat{I}_S\Transpose \mat{B} \mat{I}_S
    \in
    \RealMatrices{(\tau + 1)}{(\tau + 1)}
  $,
  and $\VolumeSampling_{\tau + 1}(\mat{B})$ is the volume sampling distribution
  prescribing to pick~$S$ with probability $\propto \det(\mat{B}_{S \times S})$.
\end{theorem}

The idea of applying volume sampling in Optimization was first proposed
in~\cite{rodomanov2020randomized} for accelerating coordinate descent methods.
It was shown that using this particular nonuniform sampling of coordinates
leads to a provable acceleration by a factor whose magnitude depends on gaps
in the spectrum of the curvature matrix.

Thus, we can interpret our basic Gradient Method (\cref{alg:GM}) with a fixed
Symmetric Polynomial Preconditioner~$\mat{P}_\tau$ as a deterministic
counterpart of the randomized coordinate descent method
from~\cite{rodomanov2020randomized} with $(\tau + 1)$-element
volume sampling of coordinates.
Correspondingly, both methods have very similar convergence properties and
theoretical efficiency estimates.

Nevertheless, this work offers several significant advantages
over~\cite{rodomanov2020randomized}.
First, in addition to the basic method, we have an accelerated one
(\cref{alg:FGM}), while the accelerated version of coordinate descent
with volume sampling is an open question.
Second, volume sampling is an expensive operation which is difficult to
carry out already when $\tau = 2$.
In contrast, the corresponding preconditioner~$\mat{P}_2$ for our gradient
methods is still computationally efficient
(see \cref{sec:DefinitionOfPreconditioners}).
Finally, as we will show next, the basic Gradient Method can be improved to
\emph{automatically} choose the best possible polynomial
preconditioner of degree~$\tau$ (including the one we have been discussing in
this section), and the resulting algorithm can easily handle much bigger values
of~$\tau$.

\newpage

\section{Krylov Subspace Preconditioning}
\label{SectionKrylov}

Our new symmetric polynomial preconditioners,
introduced in the previous section,
can be viewed as a certain family of polynomials that
we apply to our curvature matrix $\mat{B}$.
Thus, for a fixed degree $\tau > 0$, we use the matrix
$\mat{P}  =  p_{\tau}(\mat{B})$ as a preconditioner,
where  $p_{\tau}$
is a specifically constructed polynomial of degree $\tau$ such that $\mat{P} \succ 0$.

A natural question is
\textit{how optimal is this choice of a polynomial?}
Indeed, the problem of polynomial approximation has a long and rich history
with an affirmative answer provided
by the classical Chebyshev polynomials \cite{mason2002chebyshev}
for the uniform approximation bound.
We present a new adaptive algorithm that automatically
achieves the \textit{best} polynomial preconditioning. 
Then, we study what are the complexity guarantees that we can get
with this optimal approach.
In this section, we focus on the non-composite case only, i.e. the problem of unconstrained minimization
of a smooth function: $\min_{\xx \in \R^n} f(\xx)$.

\subsection{Gradient Method with Krylov Preconditioning}
\label{SubsectionKrylov}

We denote
$\mat{P}_{\aa} \Def a_0\mat{I} + a_1 \mat{B} + \ldots + a_{\tau} \mat{B}^{\tau}$,
where vector $\aa = (a_0, \ldots, a_{\tau}) \in \R^{\tau + 1}$
is a parameter. In each iteration, 
it is found by solving the linear system:
\beq \label{KrylovLinSystem}
\ba{rcl}
\aa & = & \mat{A}_{\tau}^{-1} \gg_{\tau} \;\; \in \;\; \R^{\tau + 1},
\ea
\eeq
where $\mat{A}_{\tau} = \mat{A}_{\tau}(\xx) \in \R^{(\tau + 1) \times (\tau + 1)}$
is the Gram matrix with the following structure ($0 \leq i, j \leq \tau$):
\beq \label{ADef}
\ba{rcl}
\bigl[ \mat{A}_{\tau}(\xx) \bigr]^{(i, j)} & \!\! \Def \!\! &
L \cdot \la \nabla f(\xx), \mat{B}^{i + j + 1} \nabla f(\xx) \ra,
\ea
\eeq
and $\gg_{\tau} = \gg_{\tau}(\xx) \in \R^{\tau + 1}$ is defined by
($0 \leq i \leq \tau$):
\beq \label{GDef}
\ba{rcl}
\bigl[ \gg_{\tau}(\xx) \bigr]^{(i)} & \Def & \la \nabla f(\xx), \mat{B}^{i} \nabla f(\xx) \ra.
\ea
\eeq
Note that this operation is exactly the projection 
of the direction $\frac{1}{L}\mat{B}^{-1} \nabla f(\xx_k)$ onto the \textit{Krylov subspace}
\eqref{Krylov}:
$$
\ba{rcl}
\xx_{k + 1} - \xx_k &  :=  & \argmin\limits_{\hh \in \mathcal{K}_{\xx_k}^{\tau}}
\| \hh + \frac{1}{L}\mat{B}^{-1} \nabla f(\xx_k)  \|_{\mat{B}}^2.
\ea
$$
Fortunately, for computing this projection we indeed do not need 
to invert the curvature matrix $\mat{B}$, but to solve only 
a small linear system \eqref{KrylovLinSystem} of size $\tau + 1$.
We are ready to formulate our new adaptive method.

\begin{algorithm}[H]
  \caption{Gradient Method with Krylov Preconditioning}
  \label{alg:KrylovGM}
  \begin{algorithmic}
    \STATE {\bfseries Initialization:} $\xx_0 \in \R^n$, $\tau \geq 0$, $L > 0$.
    \FOR{$k = 0, 1, \dots$}
    \STATE Form matrix $\mat{A}_{\tau}(\xx_k)$
          and vector $\gg_{\tau}(\xx_k)$ by \eqref{ADef}, \eqref{GDef}.

    \STATE Compute $\aa_{k}  =
    \mat{A}_{\tau}(\xx_k)^{-1} \gg_{\tau}(\xx_k)  \; \in \; \R^{\tau + 1}$.
    \STATE Set $\xx_{k + 1} = \xx_k - \mat{P}_{\aa_k} \nabla f(\xx_k)$.
    \ENDFOR
  \end{algorithmic}
\end{algorithm}

We prove the following optimality result.

\begin{framedtheorem} \label{TheoremKrylovGM}
  Let $\mat{P} \succ 0$ be 
  \underline{any preconditioner} that is a polynomial of degree $\tau$ of 	
 the curvature matrix:
 $$
 \ba{rcl}
 \mat{P} & = & p_{\tau}(\mat{B}), 
 \quad p_{\tau} \; \in \; \R[s],  \quad \deg(p_{\tau}) = \tau.
 \ea
 $$
 Then, for the iteration of Algorithm \ref{alg:KrylovGM} we have the global rates
  \eqref{Alg1Convex},\eqref{Alg2Convex}
  with the condition number that is attributed to $\mat{P}$ \eqref{PBound}:
  $
  \frac{\beta}{\alpha} = \frac{\beta(\mat{P})}{\alpha(\mat{P})}.
  $
\end{framedtheorem}

Hence, our method \textit{automatically} chooses the best possible 
preconditioning matrix from the polynomial class. 
Let us understand what are the bounds for $\frac{\beta}{\alpha}$
that we can achieve in this case.

\subsection{Bounds for the Condition Number}

\subsection{Proof of Proposition~\ref{PropositionCuttingEigenvalues}}

The problem of finding the best polynomial preconditioner can
be reformulated as minimizing the norm of a symmetric matrix
over the set of (positive) polynomials of a fixed degree $\tau \geq 0$:
\vspace{-1mm}
$$
\ba{rcl}
\min\limits_{p_{\tau} \in \mathcal{P}_{\tau}}
\Bigl\{ \, \gamma(p_{\tau}) & \Def & \| \mat{B} p_{\tau}(\mat{B}) - \mat{I} \| \, \Bigr\},
\ea
$$
where
$
\mathcal{P}_{\tau}  \Def 
\bigl\{
p_{\tau} \in \R[s] \, : \, \deg(p_{\tau}) = \tau, \; p_{\tau}(\mat{B}) \succ 0
\bigr\}.
$
Here we use the spectral norm to measure the size of a symmetric matrix,
and the objective can be rewritten as
\beq \label{GammaPDef}
\ba{rcl}
\gamma(p_{\tau}) & = & \max\limits_{s \in \Spec(\mat{B})} | s p_{\tau}(s) - 1 |,
\ea
\eeq
where $\Spec(\mat{B})$ is the discrete set of eigenvalues of the curvature matrix.
For any value of $\gamma := \gamma(p_{\tau})$,
our original approximation guarantee \eqref{PBound}
clearly satisfied with $\beta = 1 + \gamma$ and $\alpha = 1 - \gamma$,
and the condition number becomes\footnote{We are interested in $\gamma < 1$, since $\gamma = 1$
  trivially holds for zero polynomial.}
\vspace{-1mm}
\beq \label{CondNumberGamma}
\ba{rcl}
\frac{\beta}{\alpha} & = & \frac{1 + \gamma}{1 - \gamma}.
\ea
\eeq

Now, we take
\beq \label{AppendixQDef}
\ba{rcl}
q_{\tau}(s) & := &
\bigl( 1 - \frac{s}{\lambda_1} \bigr)
\bigl( 1 - \frac{s}{\lambda_2} \bigr)
\cdot \ldots \cdot
\bigl( 1 - \frac{s}{\lambda_{\tau}} \bigr).
\ea
\eeq

	First, note that $q_{\tau}(0) = 1$ and thus the polynomial $1 + q_{\tau}(s) \cdot (\alpha s - 1)$ is divisible by $s$.
	Hence the degree of the polynomial
	$$
	\ba{rcl}
	p_{\tau}(s) & := &  \frac{1 + q_{\tau}(s) \cdot (\alpha s - 1) }{s}
	\ea
	$$
	is exactly $\tau$.
	Then, we obtain
	\beq \label{GammaCut}
	\ba{rcl}
	\gamma & = & \max\limits_{s \in \Spec(\mat{B})} |s p_{\tau}(s) - 1 |
	\;\; = \;\;
	\max\limits_{s \in \Spec(\mat{B})} |q_{\tau}(s) \cdot (\alpha s - 1) | \\
	\\
	& \leq &
	\max\limits_{s \in \{ \lambda_{\tau + 1}, \ldots, \lambda_{n} \}} |\alpha s - 1|
	\;\;  = \;\;
	\frac{\lambda_{\tau + 1} - \lambda_n}{\lambda_{\tau + 1} + \lambda_n},
	\ea
	\eeq
	and the optimal value is $\alpha = \frac{2}{\lambda_{\tau + 1} + \lambda_n}$,
	where we  put formally $\lambda_{n + 1} \equiv \lambda_n$.
	It remains to substitute this bound into 
	$$
	\ba{rcl}
	\frac{\beta}{\alpha} & = & \frac{1 + \gamma}{1 - \gamma},
	\ea
	$$
	which is monotone in $\gamma$.
	\qed



The worst case instance for the cutting strategy
is when all the eigenvalues
except one share the same value:
$\lambda_1 = \lambda_2 = \ldots = \lambda_{n - 1} > \lambda_n$.
Indeed, then the condition number remains the same:
$
\frac{\beta}{\alpha} \, = \,
\frac{\lambda_{\tau + 1}}{\lambda_n}
\, \equiv \, \frac{\lambda_1}{\lambda_n},
$
while $\tau < n - 1$.
A better approach in such a situation would be to find a bound
from the \textit{uniform} polynomial approximation for the whole interval
$[\lambda_n, \lambda_1]$,
which is achieved with the Chebyshev polynomials \cite{polyak1987introduction,nemirovski1995information}.
Then, we can decrease the condition number
exponentially by increasing the degree $\tau$.

\begin{proposition} \label{PropositionChebyshev}
  For a fixed $0 < \epsilon < 1$, let
	$
  \tau := \bigl\lfloor \sqrt{\frac{\lambda_1}{\lambda_n}} \ln \frac{8}{\epsilon} \bigr\rfloor.
  $
  Then, taking
  $
  \mat{P}  =  p_{\tau}(\mat{B}),
  $
  where $p_{\tau}(s) := \frac{1 - Q_{\tau}(s)}{s}$
  with $Q_{\tau}$ is a normalized Chebyshev polynomial\protect\footnotemark 
  of the first kind of degree $\tau + 1$,
  the condition number is bounded by
  $
  \frac{\beta}{\alpha} \leq   1 + \epsilon.
  $
\end{proposition}

\footnotetext{See Appendix~\ref{SubsectionAppendixChebyshev} for the precise definition.}

Clearly, we can combine this technique
with the cutting strategy,
getting the best of two guarantees.

\subsection{Discussion}
\label{SubsectionDiscussion}

We see that in the case of unconstrained smooth minimization,
it is possible to achieve the guarantee of the \textit{best polynomial} of a fixed
degree $\tau$, by computing a certain projection onto the 
corresponding Krylov subspace. Namely, we can achieve 
$\frac{\beta}{\alpha} \leq \frac{\lambda_{\tau + 1}}{\lambda_n}$
(Proposition~\ref{PropositionCuttingEigenvalues}), which cuts off the top $\tau$
eigenvalues of the spectrum completely, if they are separated from the others.
At the same time, by using the Chebyshev polynomials,
we can contract a part of the spectrum \textit{uniformly}, with an 
appropriate degree $\tau$ (Proposition~\ref{PropositionChebyshev}).
It remains to be an open question where we can incorporate adaptive Krylov
preconditioning into the Fast Gradient Method, which would give us 
a further improvement
of the condition number.

Note that the linear Conjugate Gradient Method (CGM)
as applied to a convex quadratic function (thus $L = \mu = 1$), has the same guarantee
as in Theorem~\ref{TheoremKrylovGM}
with $\tau := k$, where $k$ is the iteration number.
I.e., in each iteration $k$, the method automatically achieves 
the best polynomial approximation of degree $k$.
As compared to CGM,
it remains to be the main advantage of Algorithm~\ref{alg:KrylovGM}, that
it is applicable to a wider class of nonlinear problems.

Finally, for the methods with Symmetric Polynomial Preconditioning, we obtained
the bound: $\frac{\beta}{\alpha} \leq \frac{\lambda_1}{\lambda_n} \cdot \xi_{\tau}(\blambda)$
(Theorem~\ref{TheoremPApprox}), where $\xi_{\tau}(\blambda) \leq 1$,
and $\xi_{\tau}(\blambda)  \to 0$
in the case of large spectral gap: $\lambda_1 / \lambda_{\tau + 1} \to \infty$.
This guarantee is similar to that for the cutting guarantee of Algorithm~\ref{alg:KrylovGM}.
Moreover, we can apply our family of preconditioners 
for the more general \textit{composite} optimization problems.
We also can incorporate the Symmetric Polynomial Preconditioning
into the Fast Gradient Method (Algorithm~\ref{alg:FGM}),
which significantly improves the product of the both condition numbers: 
$\frac{L}{\mu} \cdot \frac{\beta}{\alpha}$
by taking the square root (Theorem~\ref{TheoremFGM}).
Therefore, in particular practical scenarios, any of these
strategies can be more preferable.

%

 \section{Experiments}
\label{SectionExperiments}

\paragraph{Huber Loss. }

\begin{figure}[h!]
	\centering
	\includegraphics[width=0.24\linewidth]{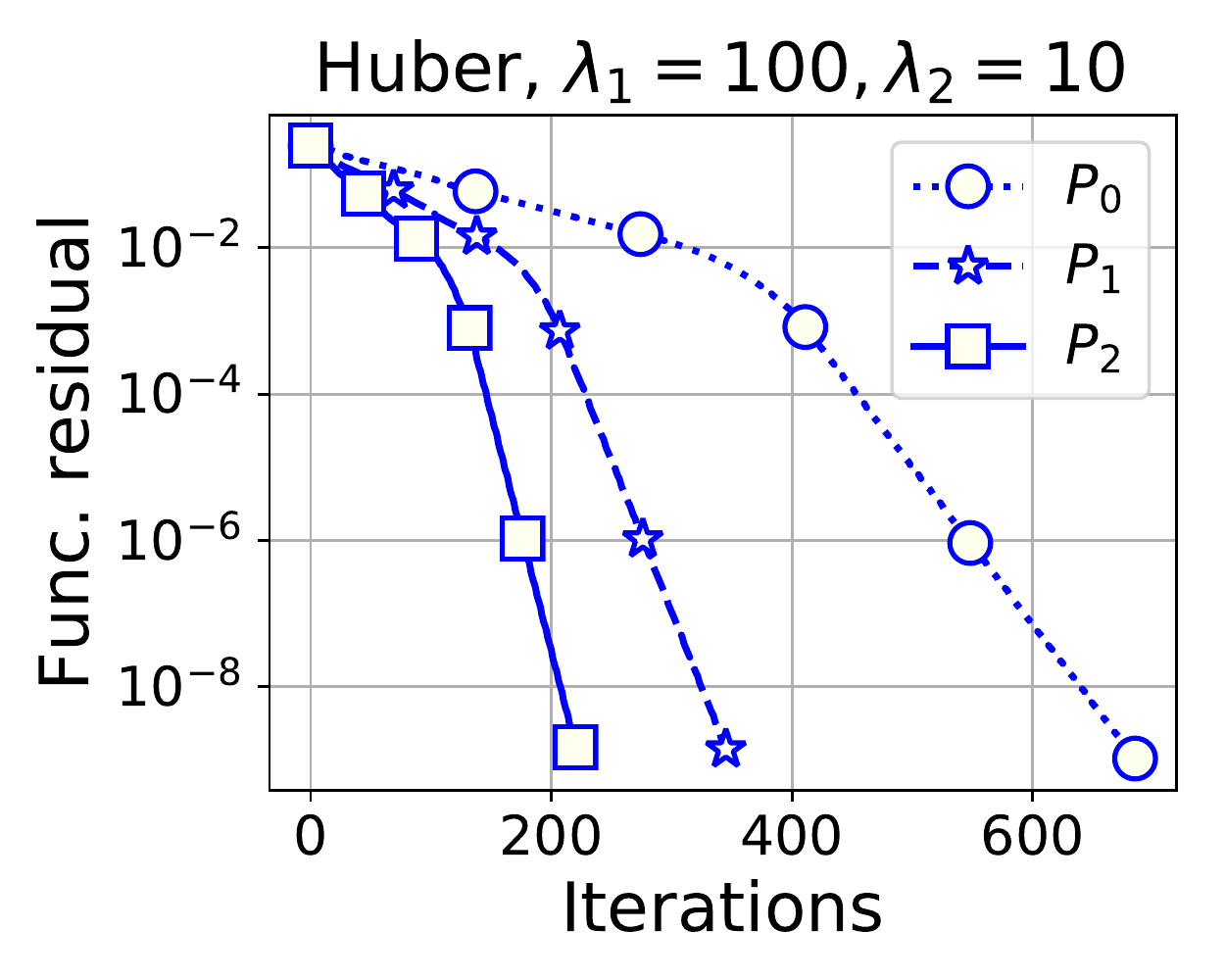}
	\includegraphics[width=0.24\linewidth]{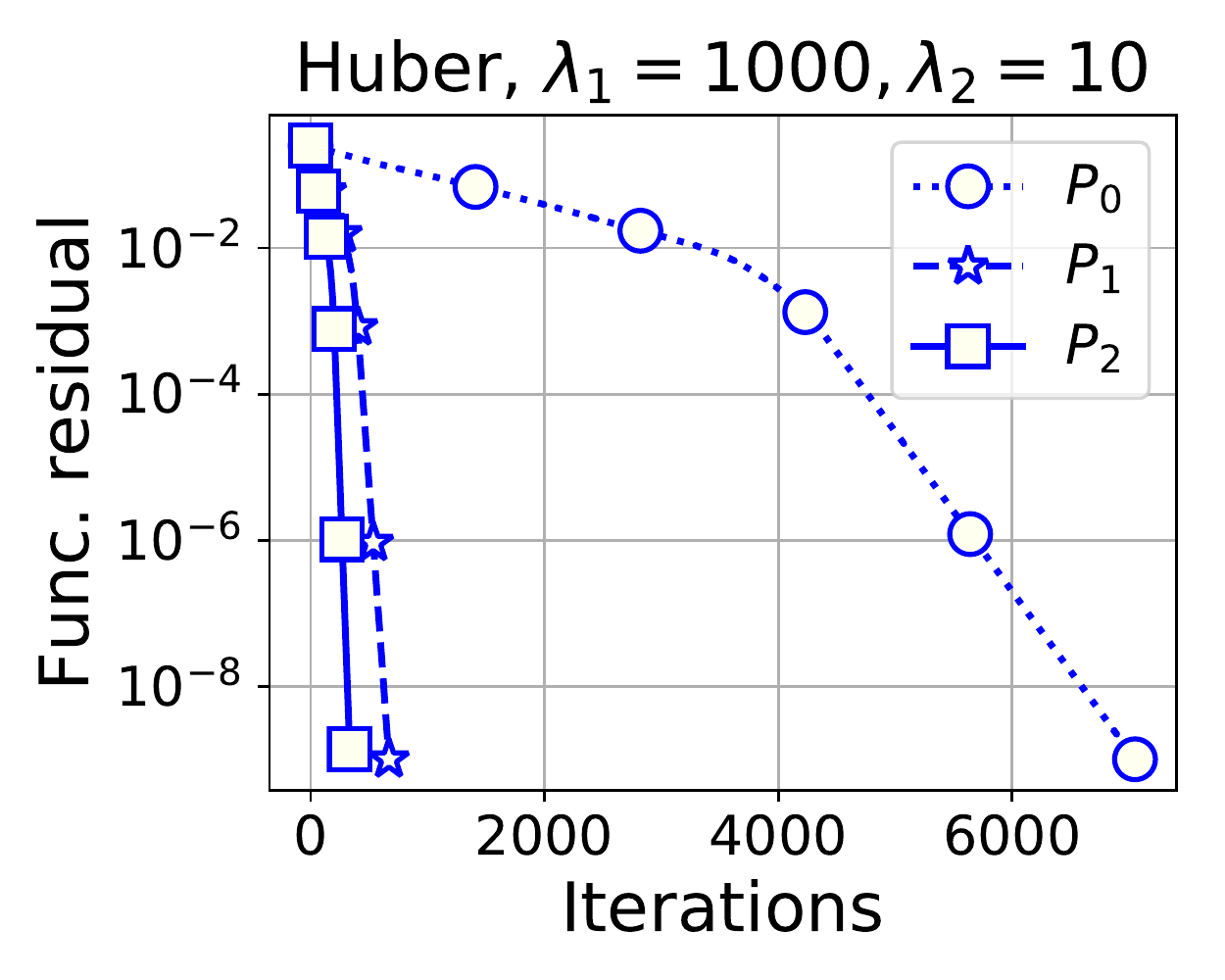}
	\includegraphics[width=0.24\linewidth]{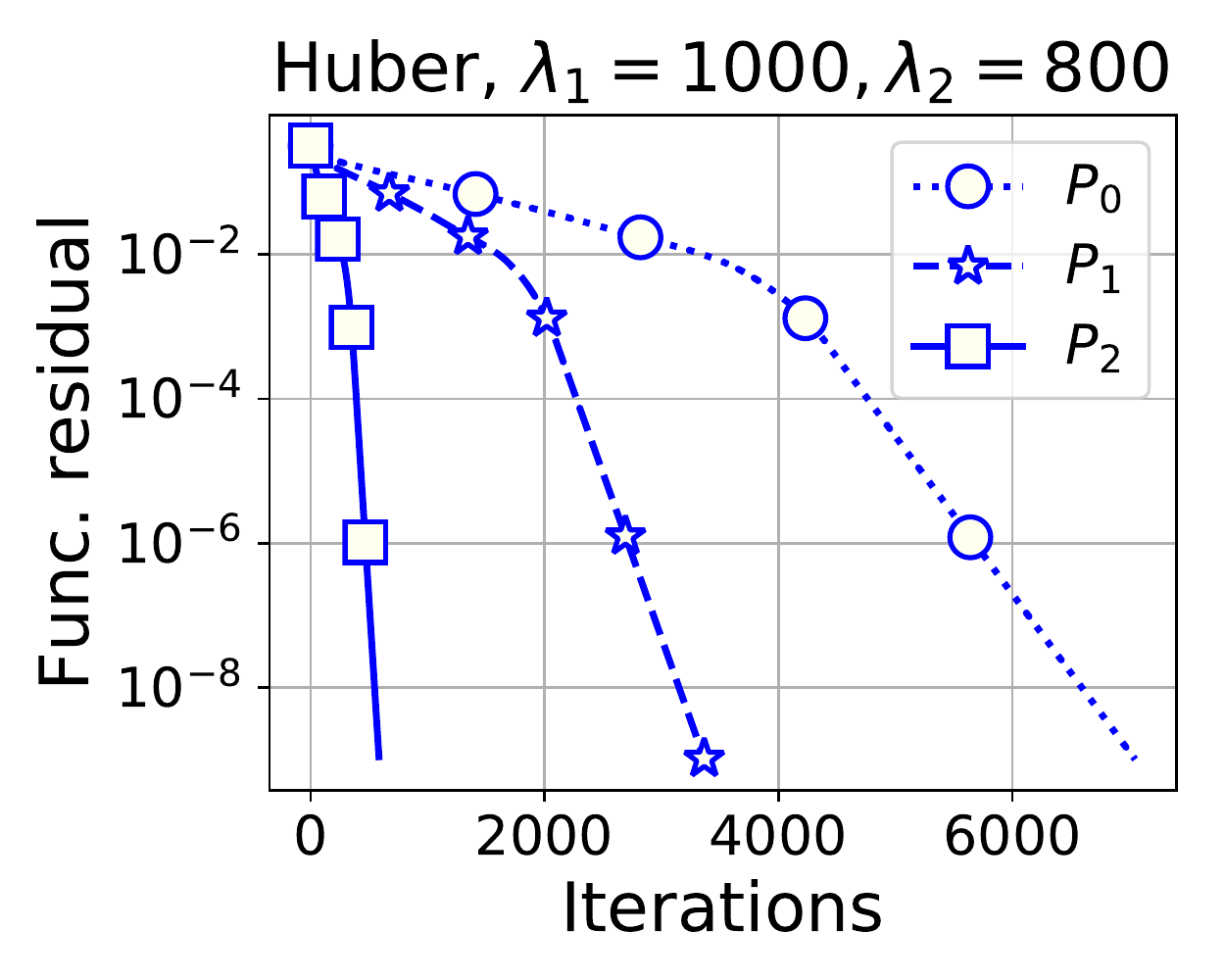}
	\includegraphics[width=0.24\linewidth]{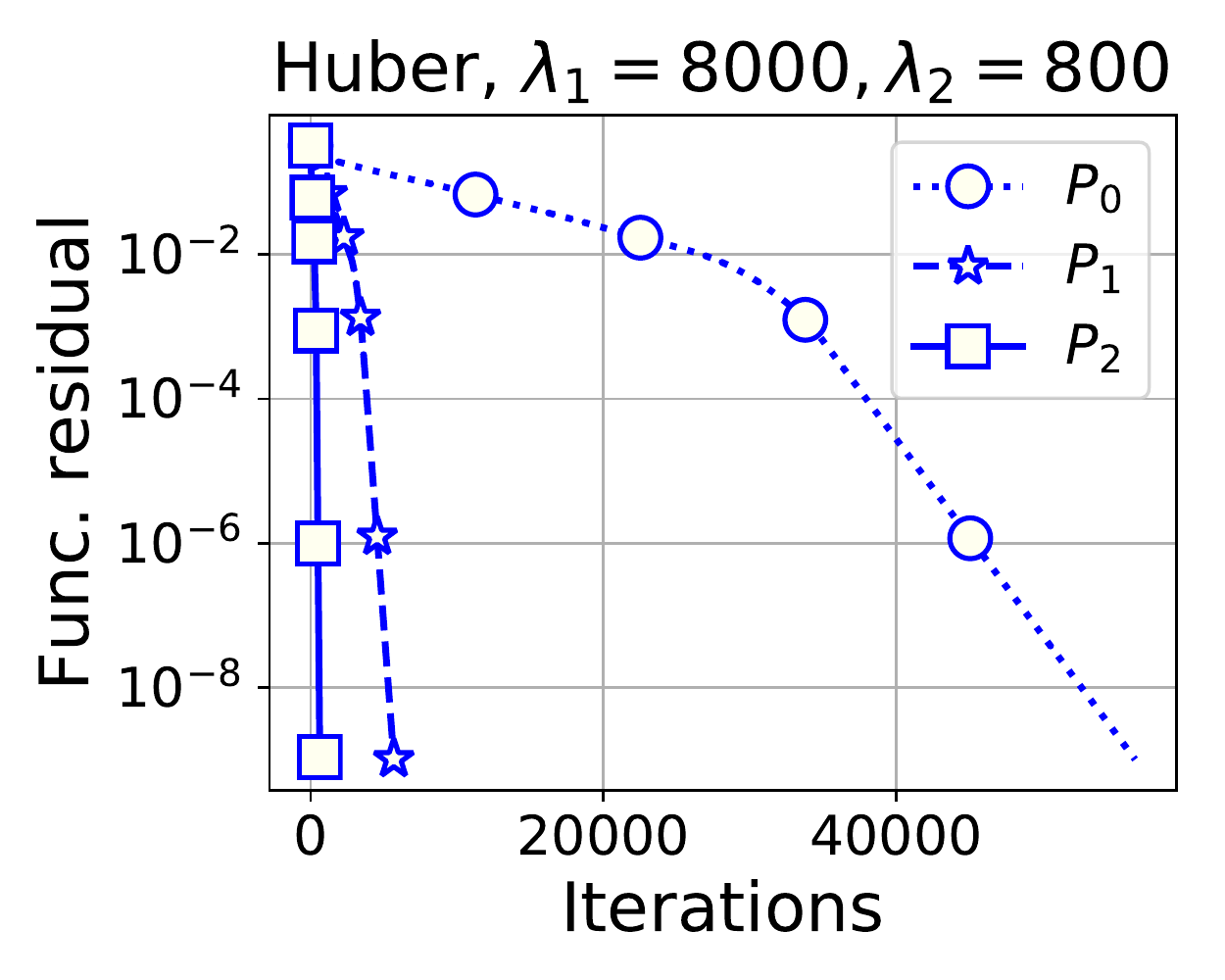}
	\caption{\small Minimizing the Huber loss by Algorithm~\ref{alg:GM} with Symmetric Polynomial Preconditioning \eqref{PRecDef}.
	$\mat{P}_0$ corresponds to the classical gradient descent without preconditioning.}
	\label{fig:Huber}
\end{figure}

Let us present an illustrative experiment,
with the regression model (Example~\ref{ExampleSeparable})
with the Huber loss function:
$$
\ba{rcl}
\phi(t) & := & \begin{cases}
\frac{t^2}{2\mu}, & \quad \text{if} \quad |t| \leq \mu, \\
|t| - \frac{\mu}{2}, & \quad \text{otherwise},
\end{cases}
\ea
$$
where $\mu := 0.1$ is a parameter.
The data is generated with a
fixed distribution of eigenvalues:
$\lambda_1 > \lambda_2 > \lambda_3 = \ldots \lambda_n = 1$.
Thus, we have two gaps between the leading eigenvalues.
We use the Gradient Method (Algorithm~\ref{alg:GM}), with
the adaptive search
to fit the parameter $M$.  
The results are shown in Figure~\ref{fig:Huber}.
Using the preconditioner $\mat{P}_1$ helps the method
to deal with the large gap between $\lambda_1$ and $\lambda_2$, while
$\mat{P}_2$ makes the method to be insensitive to the gap between $\lambda_1$ and $\lambda_3$,
as predicted by our theory. 
For example, increasing the ratio $\lambda_1 / \lambda_2$
by $10$, the performance of the methods with $\mat{P}_1$ and $\mat{P}_2$ 
becomes \textit{ten times faster} than that of the classical gradient descent.

\newpage

\paragraph{Logistic Regression.}

We examine the training of logistic regression on real data.
That problem corresponds 
to Example~\ref{ExampleSeparable} with the loss function $\phi(t) = \log(1 + e^t)$.

\begin{figure}[h!]
	\centering
	\includegraphics[width=0.70\linewidth]{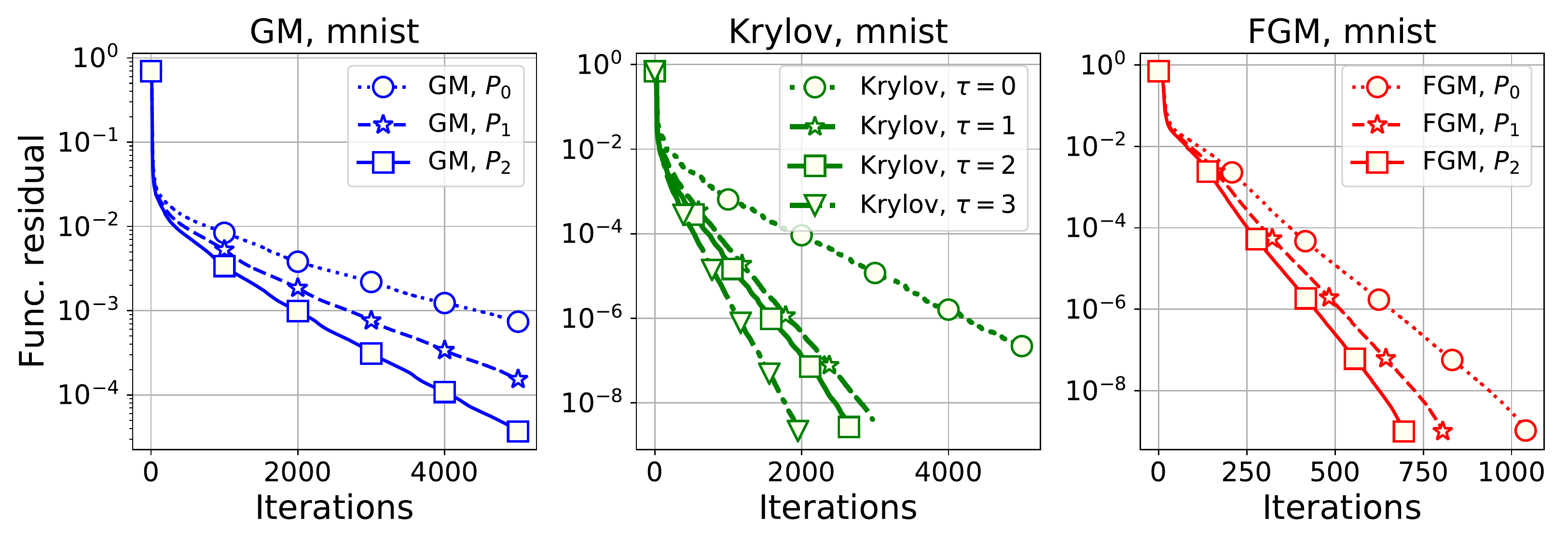}\\
	\includegraphics[width=0.70\linewidth]{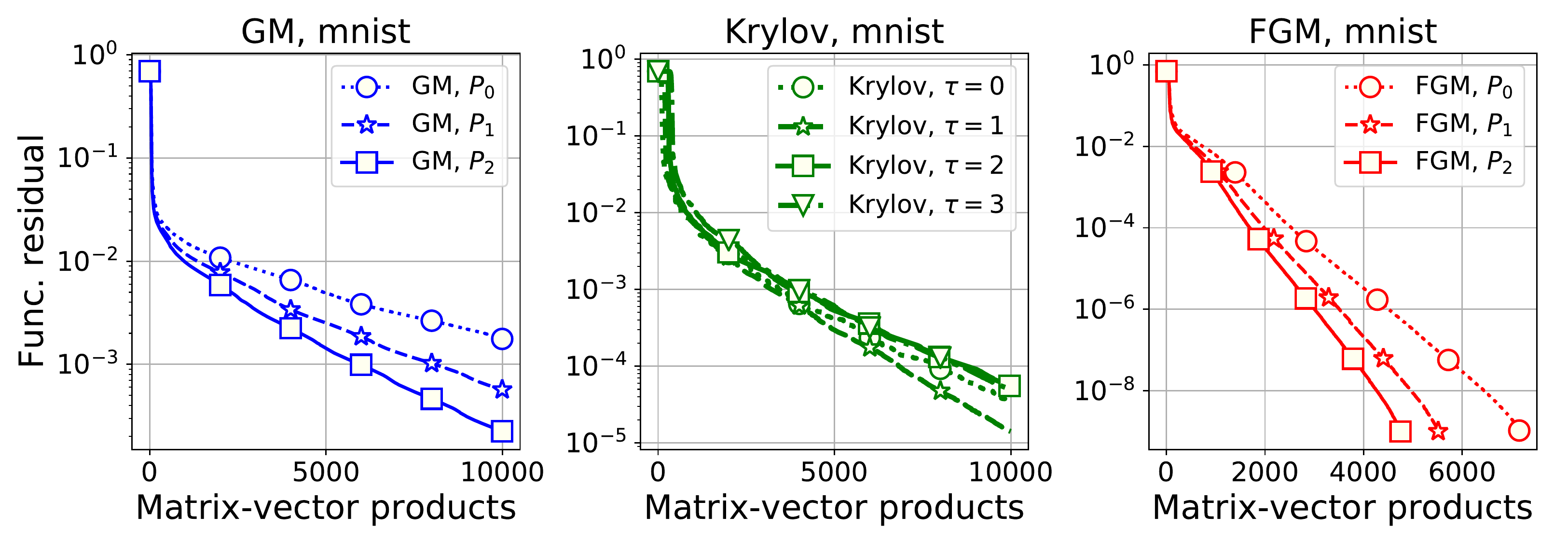}
	\caption{\small Training logistic regression with Algorithm~\ref{alg:GM} (GM)
	and
	Algorithm~\ref{alg:FGM} (FGM) employing Symmetric Polynomial Preconditioning \eqref{PRecDef};
	and with Algorithm~\ref{alg:KrylovGM} (Krylov).}
	\label{fig:Logreg}
\end{figure}

In Figure~\ref{fig:Logreg}, we see that the best convergence is achieved by the Fast Gradient Method
(FGM, Algorithm~\ref{alg:FGM}) with $\mat{P}_2$.
Using Symmetric Polynomial Preconditioning makes the methods
to converge much better (\textit{two times faster} for GM using $\mat{P}_2$ instead of $\mat{P}_0 \equiv \mat{I}$,
and about \textit{$1.5$ times faster} for FGM).
Among the versions of GM,
the most encouraging performance belongs to the Krylov preconditioning,
which is consistent with the theory.
For all the methods tested, the arithmetic cost of every iteration
remains to be at the same level.
See also Appendix~\ref{SectionExtraExperiments} for the extra experiments.

  \printbibliography

  \appendix

\newpage
\section{Extra Experiments}
\label{SectionExtraExperiments}



\paragraph{Logistic Regression. }

\begin{figure}[h!]
	\centering
	\includegraphics[width=0.49\linewidth]{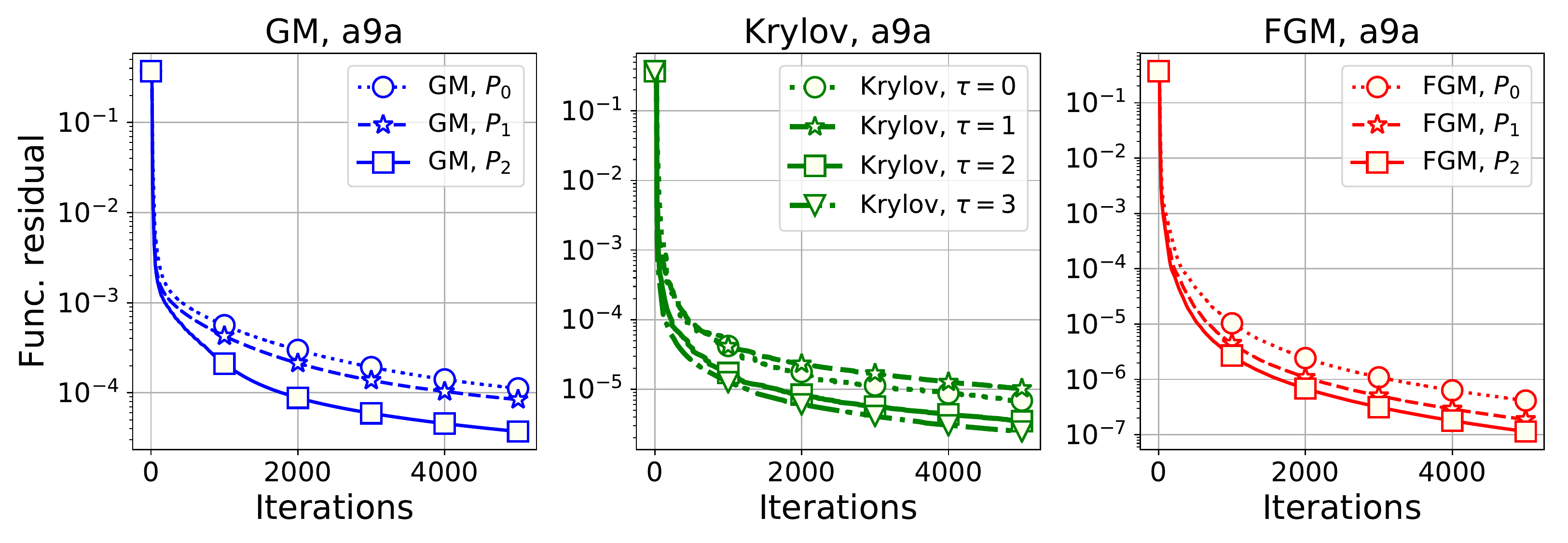}
	\includegraphics[width=0.49\linewidth]{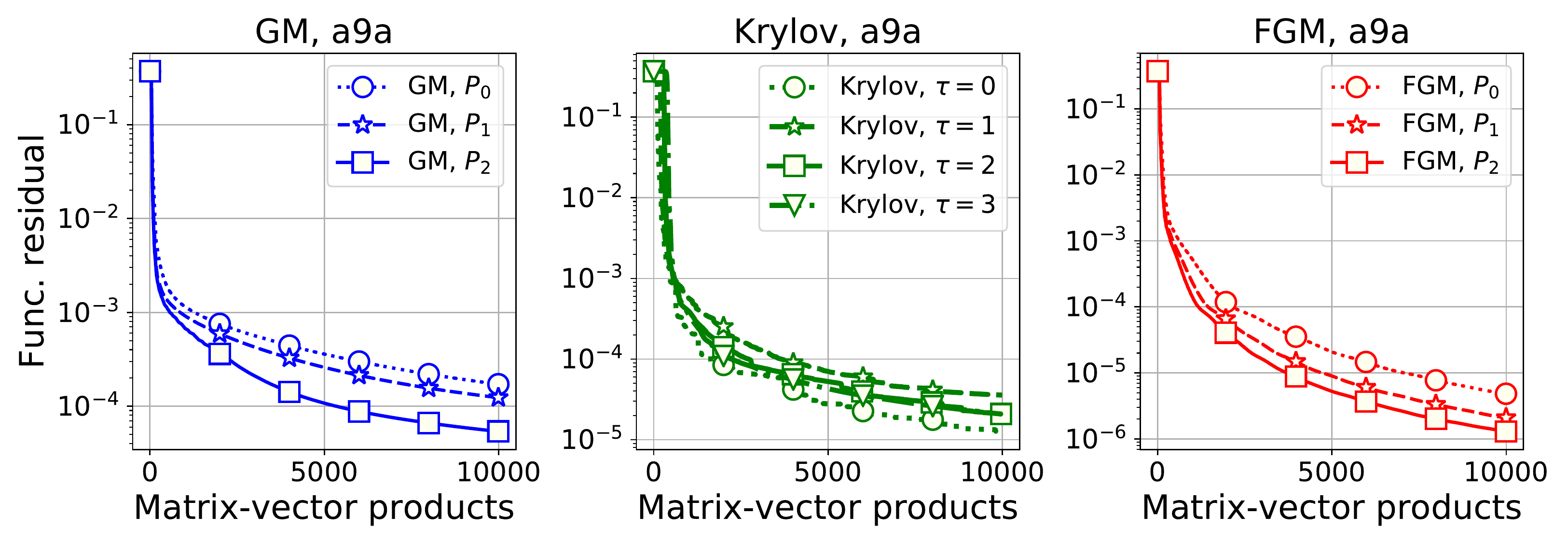}
	
	\includegraphics[width=0.49\linewidth]{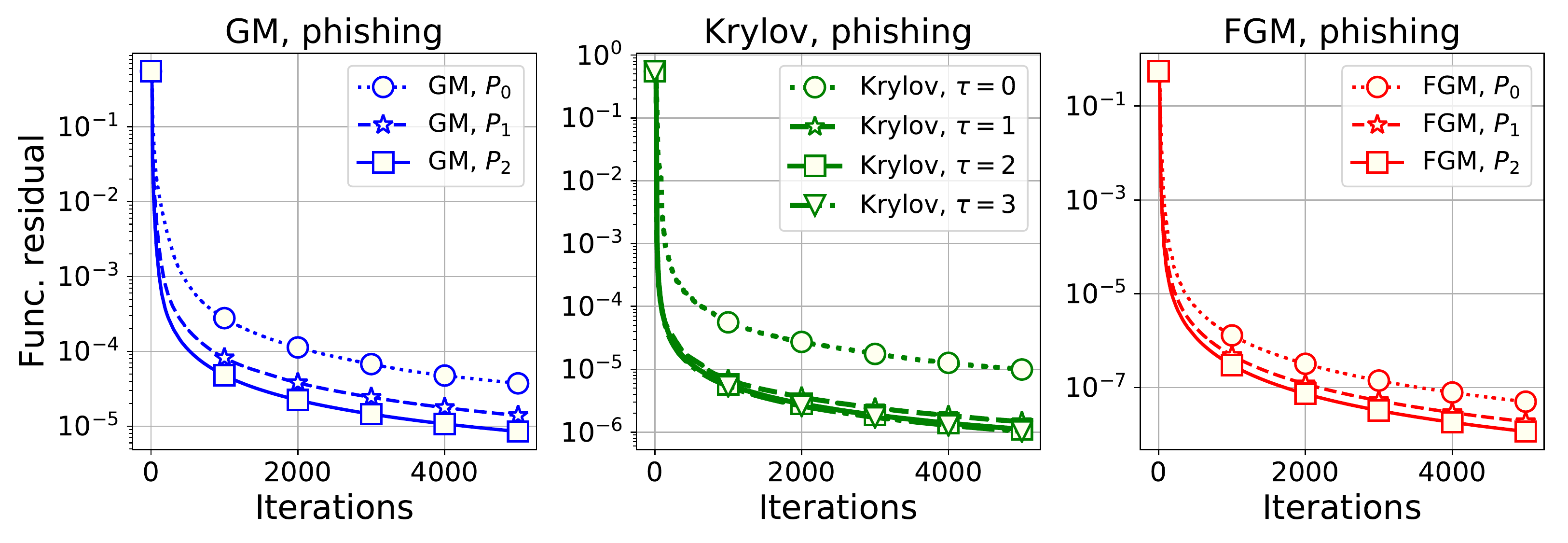}
	\includegraphics[width=0.49\linewidth]{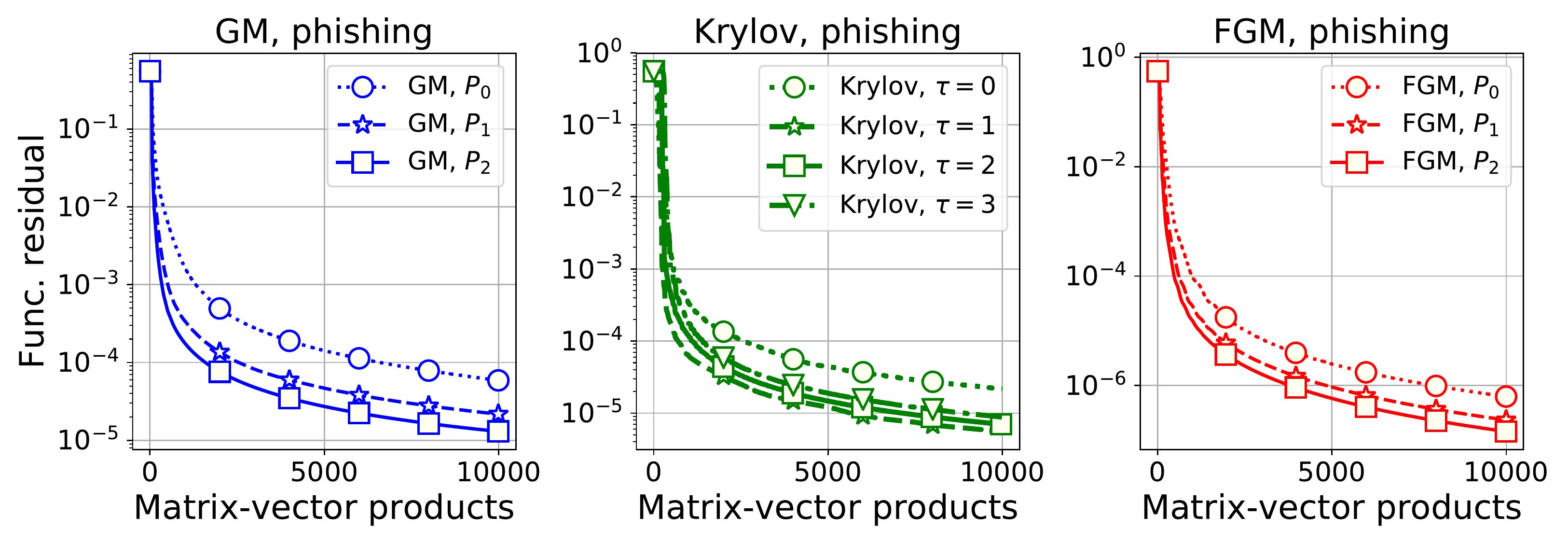}
	
	\includegraphics[width=0.49\linewidth]{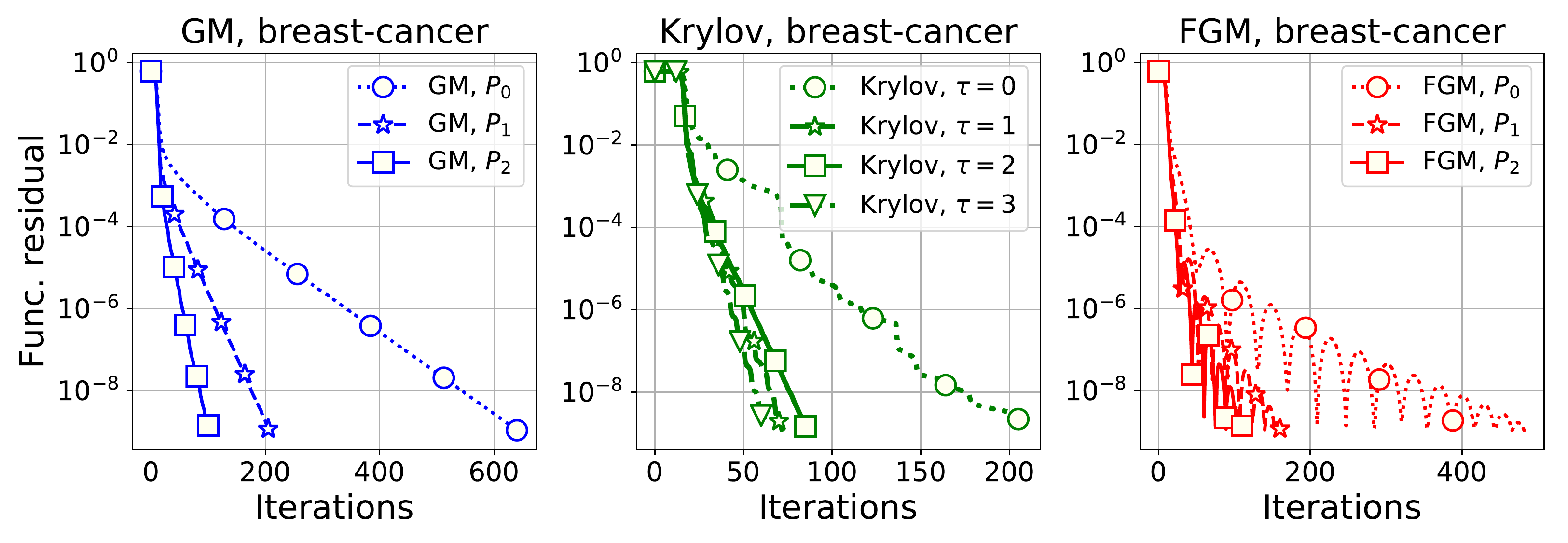}
	\includegraphics[width=0.49\linewidth]{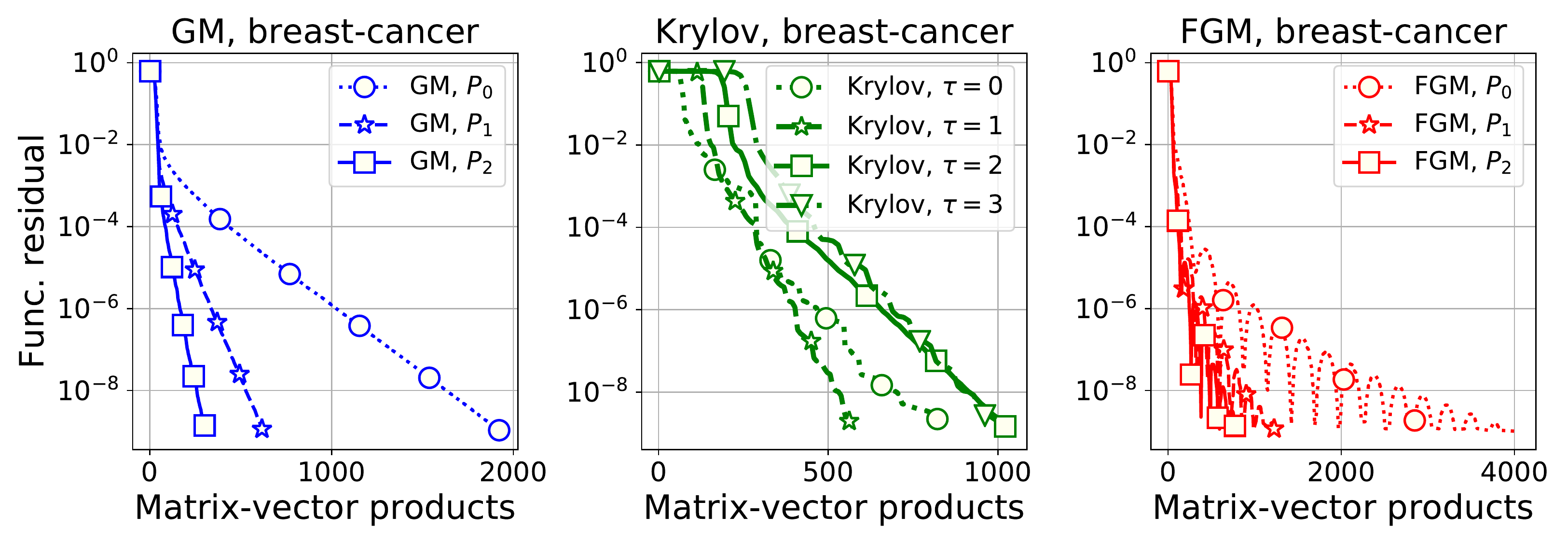}
	
	\includegraphics[width=0.49\linewidth]{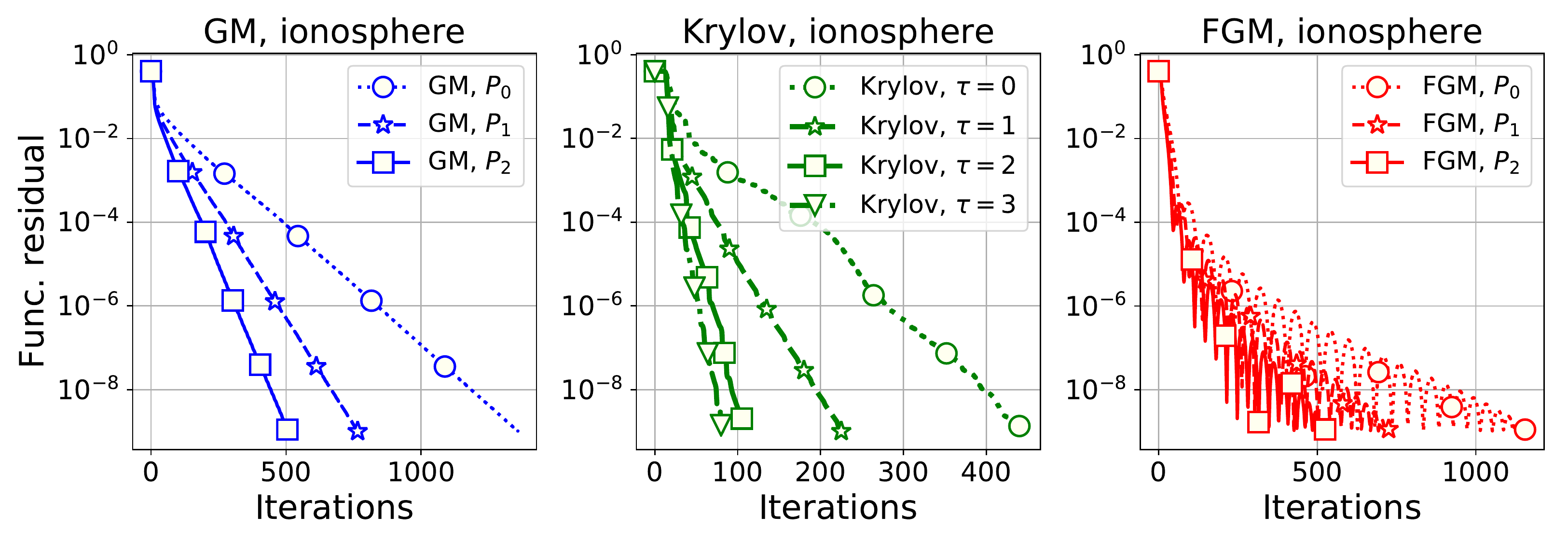}
	\includegraphics[width=0.49\linewidth]{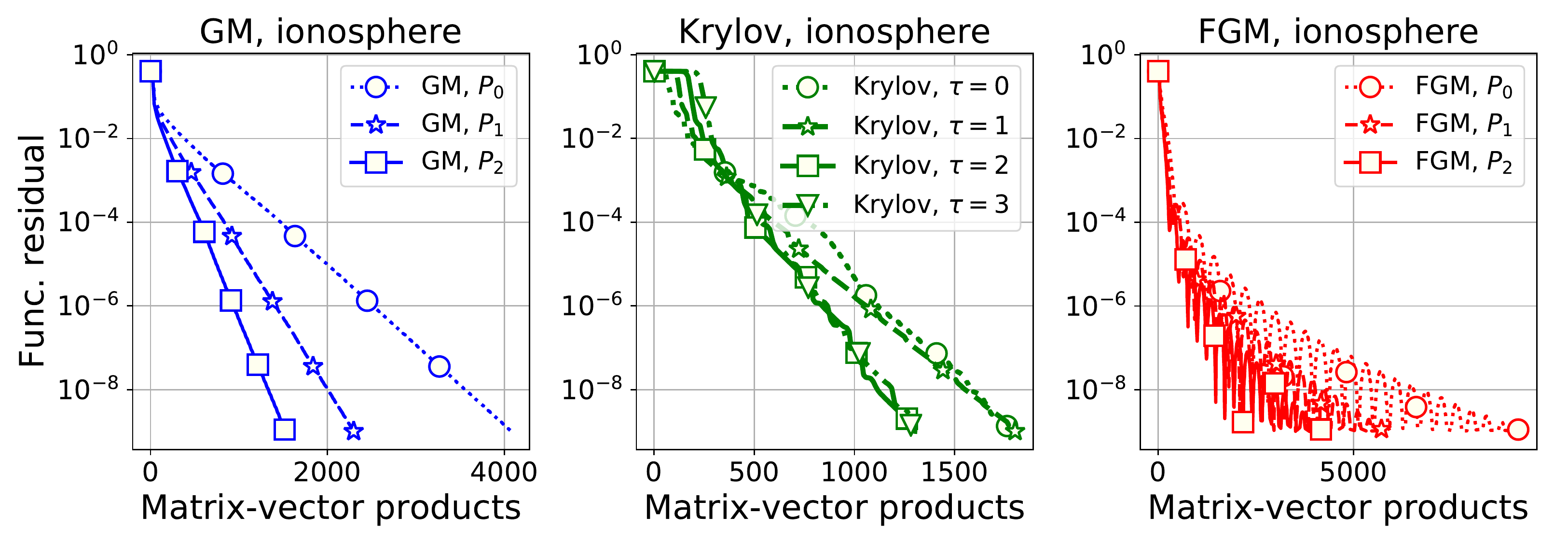}
	
		\includegraphics[width=0.49\linewidth]{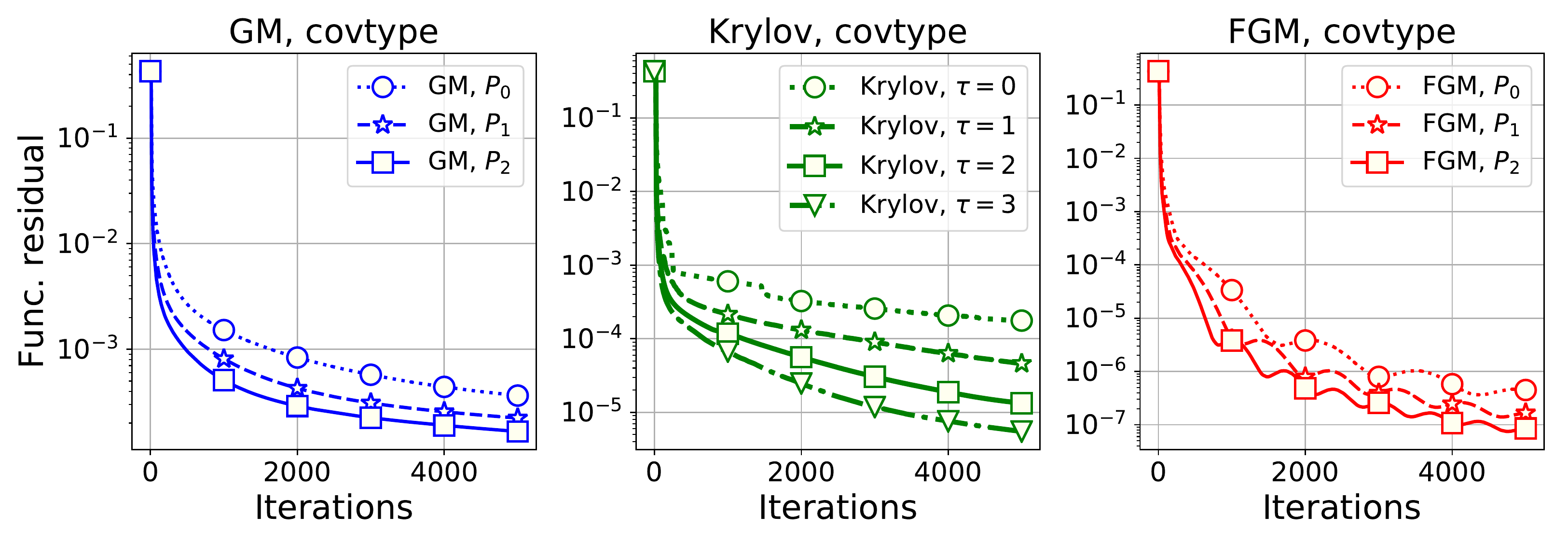}
	\includegraphics[width=0.49\linewidth]{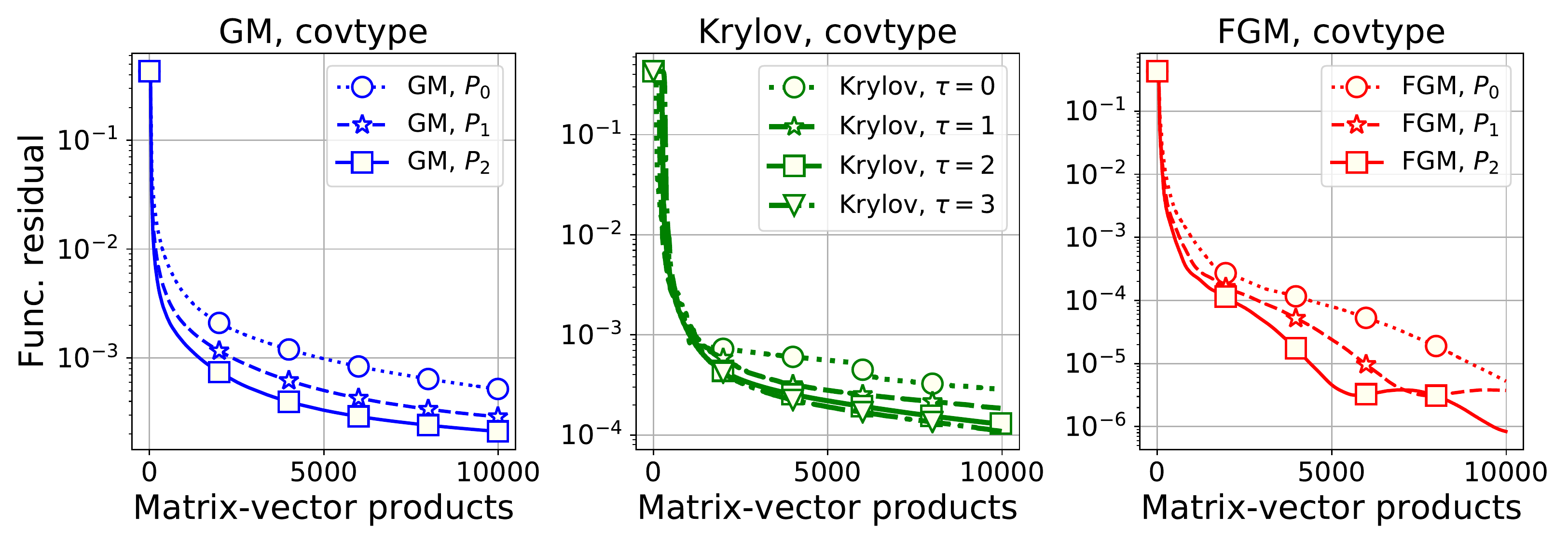}
	
	\caption{\small Training logistic regression with Algorithm~\ref{alg:GM} (GM)
		and
		Algorithm~\ref{alg:FGM} (FGM) employing Symmetric Polynomial Preconditioning \eqref{PRecDef};
		and with Algorithm~\ref{alg:KrylovGM} (Krylov).}
	\label{fig:ExtraLogreg}
\end{figure}

Let us present experimental results for our preconditioning strategies,
for the training of Logistic Regression with several real datasets.
We investigate both the number of iterations and the number of
matrix-vector products (the most difficult operation)
required to reach a certain accuracy level in the functional residual.
The results are shown in Figure~\ref{fig:ExtraLogreg}.

We see that using Symmetric Polynomial Preconditioning ($\mat{P}_1$ and $\mat{P}_2$)
significantly accelerates both the Gradient Method (GM) and the Fast Gradient Method (GM),
without extra arithmetic efforts during each iteration.
Using the Krylov preconditioning is more costly,
while it equips GM with the best possible iteration rates.


\section{Proofs}

\subsection{Proof of Theorem~\ref{TheoremGM}}
\label{sec:ProofGM}

Let us consider one iteration of the method, for some $k \geq 0$.
By definition, $\xx_{k + 1} = \argmin\limits_{\yy \in \dom \psi} \bigl\{ \Omega_k(\yy)  \bigr\}$,
where
$$
\ba{rcl}
\Omega_k(\yy) & \Def & f(\xx_k) + \la \nabla f(\xx_k), \yy - \xx_k \ra
+ \frac{M}{2}\| \yy - \xx_k \|_{\mat{P}^{-1}}^2 + \psi(\yy)
\ea
$$
is strongly convex with respect to $\mat{P}^{-1}$ norm with parameter $M := \beta L$.
Thus, we have, for any $\yy \in \dom \psi$:
\beq \label{Th1OneStep}
\ba{cl}
& \frac{M}{2}\| \yy - \xx_k \|_{\mat{P}^{-1}}^2 + F(\yy)
\; \geq \; \Omega_k(\yy)
\; \geq \;\;
\Omega_k(\xx_{k + 1})
+ \frac{M}{2}\| \yy - \xx_{k + 1} \|_{\mat{P}^{-1}}^2 \\
\\
& \; \geq \;
f(\xx_k) + \la \nabla f(\xx_k), \xx_{k + 1} - \xx_k \ra
+
\frac{L}{2}\|\xx_{k + 1} - \xx_{k}\|_{\mat{B}}^2
+  \psi(\xx_{k + 1})
+ \frac{M}{2}\| \yy - \xx_{k + 1} \|_{\mat{P}^{-1}}^2 \\
\\
& \; \overset{\eqref{FuncGlobalBounds}}{\geq} \;
F(\xx_{k + 1}) + \frac{M}{2}\|\yy - \xx_{k + 1}\|_{\mat{P}^{-1}}^2.
\ea
\eeq
Hence, substituting $\yy := \xx^{\star}$ (solution to the problem), we establish
the boundness for all iterates:
\beq \label{Th1Boundedness}
\ba{rcl}
\|\xx_{k + 1} - \xx^{\star} \|_{\mat{P}^{-1}}
& \leq &
\| \xx_k - \xx^{\star} \|_{\mat{P}^{-1}}.
\ea
\eeq
Further, let us take $\yy := \gamma_k \xx^{\star} + (1 - \gamma_k) \xx_k$, for some $\gamma_k \in [0, 1]$.
We obtain
\beq \label{Th1Gamma}
\ba{rcl}
F(\xx_{k + 1}) & \overset{\eqref{Th1OneStep}}{\leq} &
F(\gamma_k \xx^{\star} + (1 - \gamma_k) \xx_k) + \frac{\gamma_k^2 M}{2} \|\xx^{\star} - \xx_k \|_{\mat{P}^{-1}}^2 \\
\\
& \leq &
\gamma_k F^{\star} + (1 - \gamma_k) F(\xx_k)  +  \frac{\gamma_k^2 M}{2} \|\xx^{\star} - \xx_k \|_{\mat{P}^{-1}}^2.
\ea
\eeq
Now, setting $A_{k} \Def k \cdot (k + 1)$, $a_{k + 1} \Def A_{k + 1} - A_k = 2(k + 1)$, and $\gamma_k := \frac{a_{k + 1}}{A_{k + 1}} = \frac{2}{k + 2}$,
we obtain
\beq \label{Th1PreTelescoped}
\ba{rcl}
A_{k + 1}\bigl( F(\xx_{k + 1}) - F^{\star} \bigr)
& \overset{\eqref{Th1Gamma}}{\leq} &
A_k \bigl( F(\xx_k) - F^{\star} \bigr) + \frac{a_{k + 1}^2}{A_{k + 1}} \cdot \frac{M}{2} \|\xx^{*} - \xx_k\|_{\mat{P}^{-1}}^2 \\
\\
& \overset{\eqref{Th1Boundedness}}{\leq} &
A_k \bigl( F(\xx_k) - F^{\star} \bigr) + \frac{a_{k + 1}^2}{A_{k + 1}} \cdot \frac{M}{2} \|\xx^{*} - \xx_0\|_{\mat{P}^{-1}}^2 \\
\\
& \overset{\eqref{PBound}}{\leq} &
A_k \bigl( F(\xx_k) - F^{\star} \bigr) + \frac{a_{k + 1}^2}{A_{k + 1}} \cdot \frac{\beta}{\alpha} \cdot \frac{L}{2} \|\xx^{*} - \xx_0\|_{\mat{B}}^2.
\ea
\eeq
Telescoping this bound for the first $k$ iterations, we get
$$
\ba{rcl}
F(\xx_k) - F^{\star} & \overset{\eqref{Th1PreTelescoped}}{\leq} &
\frac{\beta}{\alpha} \cdot \frac{L}{2} \|\xx^{*} - \xx_0\|_{\mat{B}}^2
\cdot \frac{1}{A_k}
\sum\limits_{i = 1}^k \frac{a_i^2}{A_i}
\;\; = \;\;
\cO\Bigl( \frac{\beta}{\alpha} \cdot \frac{L}{2k} \|\xx^{*} - \xx_0\|_{\mat{B}}^2 \Bigr).
\ea
$$
To prove the linear rate for the strongly convex case, we continue as follows
$$
\ba{rcl}
F(\xx_{k + 1}) & \overset{\eqref{Th1Gamma},\eqref{FuncGlobalBounds}}{\leq} &
\gamma_k F^{\star} + (1 - \gamma_k) F(\xx_k)  + \gamma_k^2 \cdot \frac{\beta L }{\alpha \mu}
\cdot \bigl( F(\xx_k) - F^{\star} \bigr).
\ea
$$
Choosing $\gamma_k := \frac{\alpha \mu}{2 \beta L } < 1$, we get the exponential rate
$$
\ba{rcl}
F(\xx_{k + 1}) - F^{\star} & \leq & \Bigl(1 - \frac{\alpha \mu}{4 \beta L} \Bigr) \bigl( F(\xx_k) - F^{\star} \bigr),
\ea
$$
which completes the proof. \qed

\subsection{Proof of \cref{TheoremFGM}}

Let $\Vector{x} \in \EffectiveDomain \psi$ and $k \geq 0$ be arbitrary.
From~\eqref{FuncGlobalBounds}, \eqref{PBound}, and the fact that
$\rho = \alpha \mu$, it follows that
\[
  F(\Vector{x})
  =
  f(\Vector{x}) + \psi(\Vector{x})
  \geq
  \ell_k(\Vector{x})
  +
  \frac{\rho}{2} \RelativeNorm{\Vector{x} - \Vector{y}_k}{\mat{P}^{-1}}^2,
  \qquad
  \ell_k(\Vector{x})
  \Def
  f(\Vector{y}_k)
  +
  \InnerProduct{\Gradient f(\Vector{y}_k)}{\Vector{x} - \Vector{y}_k}
  +
  \psi(\Vector{x}).
\]
Hence,
\begin{equation}
  \label{eq:PreliminaryStep}
  \begin{aligned}
    \hspace{2em}&\hspace{-2em}
    A_k F(\Vector{x}_k)
    +
    a_{k + 1} F(\Vector{x})
    +
    \frac{1 + \rho A_k}{2}
    \RelativeNorm{\Vector{x} - \Vector{v}_k}{\mat{P}^{-1}}^2
    \\
    &\geq
    A_k \ell_k(\Vector{x}_k)
    +
    a_{k + 1} \ell_k(\Vector{x})
    +
    \frac{1 + \rho A_k}{2}
    \RelativeNorm{\Vector{x} - \Vector{v}_k}{\mat{P}^{-1}}^2
    +
    \frac{\rho a_{k + 1}}{2}
    \RelativeNorm{\Vector{x} - \Vector{y}_k}{\mat{P}^{-1}}^2
    \\
    &\geq
    A_k \ell_k(\Vector{x}_k)
    +
    a_{k + 1} \ell_k(\Vector{x})
    +
    \frac{1 + \rho A_{k + 1}}{2}
    \RelativeNorm{\Vector{x} - \hat{\Vector{v}}_k}{\mat{P}^{-1}}^2
    \Def
    \zeta_k(\Vector{x}),
  \end{aligned}
\end{equation}
where the final inequality follows from the convexity of the squared norm
and the fact that, according to our definitions,
\[
  \frac{
    (1 + \rho A_k) \Vector{v}_k + \rho a_{k + 1} \Vector{y}_k
  }{
    1 + \rho A_{k + 1}
  }
  =
  (1 - \omega_k) \Vector{v}_k + \omega_k \Vector{y}_k
  =
  (1 - \omega_k) \Vector{v}_k
  +
  \omega_k [(1 - \theta_k) \Vector{x}_k + \theta_k \hat{\Vector{v}}_k]
  =
  \hat{\Vector{v}}_k.
\]
Note that $\zeta_k$ is a $(1 + \rho A_{k + 1})$-strongly convex function
w.r.t.\ $\RelativeNorm{\cdot}{\mat{P}^{-1}}$,
and $\Vector{v}_{k + 1}$ is precisely its minimizer.
Therefore,
\begin{equation}
  \label{eq:OptimalityConditionForProxCenter}
  \zeta_k(\Vector{x})
  \geq
  \zeta_k(\Vector{v}_{k + 1})
  +
  \frac{1 + \rho A_{k + 1}}{2}
  \RelativeNorm{\Vector{x} - \Vector{v}_{k + 1}}{\mat{P}^{-1}}^2.
\end{equation}
Since $\ell_k$ is a convex function, we have, by our definition
of~$\Vector{x}_{k + 1}$,
\[
  A_k \ell_k(\Vector{x}_k) + a_{k + 1} \ell_k(\Vector{v}_{k + 1})
  \geq
  A_{k + 1} \ell_k(\Vector{x}_{k + 1}).
\]
On the other hand, by the definition of $\Vector{x}_{k + 1}$ and $\Vector{y}_k$,
\begin{align*}
  \Vector{x}_{k + 1} - \Vector{y}_k
  =
  \theta_k (\Vector{v}_{k + 1} - \hat{\Vector{v}}_k)
  =
  \frac{a_{k + 1}}{A_{k + 1}} (\Vector{v}_{k + 1} - \hat{\Vector{v}}_k).
\end{align*}
Therefore,
\begin{align*}
  \zeta_k(\Vector{v}_{k + 1})
  &=
  A_k \ell_k(\Vector{x}_k) + a_{k + 1} \ell_k(\Vector{v}_{k + 1})
  +
  \frac{1 + \rho A_{k + 1}}{2} \RelativeNorm{\Vector{v}_{k + 1} - \hat{\Vector{v}}_k}{\mat{P}^{-1}}^2
  \\
  &\geq
  A_{k + 1} \Bigl[
    \ell_k(\Vector{x}_{k + 1})
    +
    \frac{A_{k + 1} (1 + \rho A_{k + 1})}{2 a_{k + 1}^2}
    \RelativeNorm{\Vector{x}_{k + 1} - \Vector{y}_k}{\mat{P}^{-1}}^2
  \Bigr].
\end{align*}
In view of our choice of~$a_{k + 1}$, we have the following identity:
\begin{equation}
  \label{eq:EquationForCoefficient}
  \frac{M a_{k + 1}^2}{A_{k + 1}} = 1 + \rho A_{k + 1}.
\end{equation}
Combining this with the fact that $M = \beta L$ and using~\eqref{PBound}
and~\eqref{FuncGlobalBounds}, we get
\begin{align*}
  \zeta_k(\Vector{v}_{k + 1})
  &\geq
  A_{k + 1} \Bigl[
    \ell_k(\Vector{x}_{k + 1})
    +
    \frac{M}{2} \RelativeNorm{\Vector{x}_{k + 1} - \Vector{y}_k}{\mat{P}^{-1}}^2
  \Bigr]
  \geq
  A_{k + 1} \Bigl[
    \ell_k(\Vector{x}_{k + 1})
    +
    \frac{L}{2} \RelativeNorm{\Vector{x}_{k + 1} - \Vector{y}_k}{B}^2
  \Bigr]
  \\
  &=
  A_{k + 1} \Bigl[
    f(\Vector{y}_k) + \InnerProduct{\Gradient f(\Vector{y}_k)}{\Vector{x}_{k + 1} - \Vector{y}_k}
    +
    \frac{L}{2} \RelativeNorm{\Vector{x}_{k + 1} - \Vector{y}_k}{B}^2
    +
    \psi(\Vector{x}_{k + 1})
  \Bigr]
  \geq
  A_{k + 1} F(\Vector{x}_{k + 1}).
\end{align*}

Substituting the above bound
into~\eqref{eq:OptimalityConditionForProxCenter},
and that one into~\eqref{eq:OptimalityConditionForProxCenter},
we thus obtain
\[
  A_k F(\Vector{x}_k)
  +
  a_{k + 1} F(\Vector{x})
  +
  \frac{1 + \rho A_k}{2}
  \RelativeNorm{\Vector{x} - \Vector{v}_k}{\mat{P}^{-1}}^2
  \geq
  A_{k + 1} F(\Vector{x}_{k + 1})
  +
  \frac{1 + \rho A_{k + 1}}{2}
  \RelativeNorm{\Vector{x} - \Vector{v}_{k + 1}}{\mat{P}^{-1}}^2.
\]
This inequality is valid for any $k \geq 0$.

Fixing an arbitrary $k \geq 1$ and summing up the previous inequalities
for all indices $k' = 0, \ldots, k - 1$, we get
\[
  A_k F(\Vector{x}_k)
  \leq
  A_k F(\Vector{x})
  +
  \frac{1 + \rho A_0}{2}
  \RelativeNorm{\Vector{x} - \Vector{v}_0}{\mat{P}^{-1}}^2
  =
  A_k F(\Vector{x})
  +
  \frac{1}{2} \RelativeNorm{\Vector{x} - \Vector{x}_0}{\mat{P}^{-1}}^2.
\]
Substituting further $\Vector{x} = \Vector{x}^{\star}$ (an optimal solution)
and using~\eqref{PBound}, gives us the following convergence rate estimate:
\begin{equation}
  \label{eq:PreliminaryEstimate}
  F(\Vector{x}_k) - F^{\star}
  \leq
  \frac{\RelativeNorm{\Vector{x}^{\star} - \Vector{x}_0}{\mat{P}^{-1}}^2}{2 A_k}
  \leq
  \frac{\RelativeNorm{\Vector{x}^{\star} - \Vector{x}_0}{B}^2}{2 \alpha A_k}.
\end{equation}

To complete the proof, it remains to use standard lower bounds on~$A_k$
(see~\cite{nesterov2018lectures}).
Specifically, dropping the second term from the right-hand side
of~\eqref{eq:EquationForCoefficient} and rearranging, we obtain,
for any $k \geq 0$,
\[
  \sqrt{\frac{A_{k + 1}}{M}}
  \leq
  a_{k + 1}
  =
  A_{k + 1} - A_k
  =
  (\sqrt{A_{k + 1}} - \sqrt{A_k} \,) (\sqrt{A_{k + 1}} + \sqrt{A_k} \,)
  \leq
  2 (\sqrt{A_{k + 1}} - \sqrt{A_k} \,) \sqrt{A_{k + 1}}.
\]
Cancelling $\sqrt{A_{k + 1}}$ on both sides and using the fact that $A_0 = 0$,
we obtain, for any $k \geq 1$,
\[
  \sqrt{A_k} \geq \frac{k}{2 \sqrt{M}}.
\]
Squaring both sides, substituting the resulting inequality
into~\eqref{eq:PreliminaryEstimate} and replacing $M = \beta L$,
we get~\eqref{Alg2Convex}.

When $\mu > 0$, we can drop the first term from the
right-hand side of~\eqref{eq:EquationForCoefficient}.
This gives us
\[
  a_{k + 1}^2 \geq \frac{\rho}{M} A_{k + 1}^2.
\]
Hence, for any $k \geq 0$,
\[
  A_{k + 1} - A_k
  =
  a_{k + 1}
  \geq
  q A_{k + 1},
  \qquad
  q \Def \sqrt{\frac{\rho}{M}} \leq 1,
\]
or, equivalently,
\[
  A_{k + 1}
  \geq
  \frac{A_k}{1 - q}.
\]
Consequently, for any $k \geq 1$,
\[
  A_k
  \geq
  \frac{A_1}{(1 - q)^{k - 1}}
  \geq
  \frac{1}{M (1 - q)^{k - 1}},
\]
where the final inequality is due to~\eqref{eq:EquationForCoefficient}
combined with the fact that $A_0 = 0$.
Substituting this inequality into~\eqref{eq:PreliminaryEstimate}
and replacing $M = \beta L$, $\rho = \alpha \mu$, we get
the second bound from \cref{TheoremFGM}.
\qed

\subsection{Proof of Lemma~\ref{LemmaSpec}}

Let us denote by $u_k(\aa)$ the $k$-th \textit{power sum} of 
the variables:
$$
\ba{rcl}
u_k(\aa) & \Def & \sum\limits_{i = 1}^{n - 1} a_i^{k}, \qquad 
\forall \aa \in \R^{n - 1}.
\ea
$$
Then, the classical Newton-Girard identities (see, e.g. \cite{kalman2000matrix}) state
the following relation between the elementary symmetric polynomials:
\beq \label{NGIdent}
\ba{rcl}
\sigma_{\tau}(\aa)
& \equiv &
\frac{1}{k} \sum\limits_{i = 1}^\tau (-1)^{i - 1} \sigma_{\tau - i}(\aa) \cdot u_i(\aa).
\ea
\eeq
Note that for the matrix $\mat{U}_{\tau} \Def \tr(\mat{B}^{\tau}) \mat{I} - \mat{B}^{\tau}$,
the following spectral decomposition holds:
\beq \label{UtauSpec}
\ba{rcl}
\mat{U}_{\tau} & = & 
\mat{Q} \Diag\Bigl(  \sum\limits_{i = 1}^n \lambda_i^{\tau} - \lambda_1^{\tau},
\, \sum\limits_{i = 1}^n \lambda_i^{\tau}  - \lambda_2^{\tau}, \ldots, 
\sum\limits_{i = 1}^n \lambda_i^{\tau} - \lambda_n^{\tau} \Bigr) \mat{Q}^{\top} \\
\\
& = & 
\mat{Q} \Diag\Bigl(  
u_{\tau}(\blambda_{-1}), u_{\tau}(\blambda_{-2}), \ldots, u_{\tau}(\blambda_{-n})
\Bigr) \mat{Q}^{\top}. 
\ea
\eeq
Now, the identity that we need to prove is
\beq \label{SpecProof}
\ba{rcl}
\mat{P}_{\tau}
& = & 
\mat{Q} \Diag\Bigl(
\sigma_{\tau}(\blambda_{-1}), \,
\sigma_{\tau}(\blambda_{-2}), \ldots, \sigma_{\tau}(\blambda_{-n})
 \Bigr) \mat{Q}^{\top}.
\ea
\eeq
We justify \eqref{SpecProof} by induction. By definition, $\mat{P}_0 \Def \mat{I}$
and $\sigma_{0}(\aa) \equiv 1$, therefore \eqref{SpecProof} holds for $\tau = 0$, which is our base.
Let us fix $\tau \geq 1$ and assume that \eqref{SpecProof} is true for all smaller indices. Then,
$$
\ba{rcl}
\mat{P}_{\tau} & \Def &
\frac{1}{\tau} \sum\limits_{i = 1}^{\tau} (-1)^{i - 1} \mat{P}_{\tau - i} \mat{U}_{i} \\
\\
& \overset{\eqref{SpecProof},\eqref{UtauSpec}}{=} &
\mat{Q} \Diag\Bigl( \,
\sum\limits_{i = 1}^{\tau} (-1)^{i - 1} \sigma_{\tau - i}(\blambda_{-1}) \cdot u_{i}(\blambda_{-1}),
\, \ldots, \, 
\sum\limits_{i = 1}^{\tau} (-1)^{i - 1} \sigma_{\tau - i}(\blambda_{-n}) \cdot u_{i}(\blambda_{-n})
\, \Bigr) \mat{Q}^{\top} \\
\\
& \overset{\eqref{NGIdent}}{=} &
\mat{Q} \Diag\Bigl( \,
\sigma_{\tau}(\blambda_{-1}),
\, \ldots, \,
\sigma_{\tau}(\blambda_{-n})   \Bigr) \mat{Q}^{\top}.
\ea
$$
Hence, \eqref{SpecProof} is proven for all $0 \leq \tau \leq n - 1$. \qed

\subsection{Proof of Theorem~\ref{TheoremPApprox}}

By Lemma~\ref{LemmaSpec}, we have
the following representation of our preconditioner:
$$
\ba{rcl}
\mat{P}_{\tau} 
& = & \mat{Q} \Diag(
\sigma_{\tau}(\blambda_{-1}), \sigma_{\tau}(\blambda_{-2}), \ldots, 
\sigma_{\tau}(\blambda_{-n}))
 ) \mat{Q}^{\top}.
\ea
$$
It is easy to see that, for the spectrum of the matrix
$$
\ba{rcl}
\mat{B}^{1/2} \mat{P}_{\tau} \mat{B}^{1/2}
& = &
\mat{Q} \Diag\Bigl(  
\lambda_1 \cdot \sigma_{\tau} (\blambda_{-1}) ,
\lambda_2 \cdot \sigma_{\tau} (\blambda_{-2}) ,
\ldots,
\lambda_n \cdot \sigma_{\tau} (\blambda_{-n}) 
 \Bigr) \mat{Q}^{\top},
\ea
$$
it holds:
\beq \label{th4_2_bound}
\ba{rcl}
\lambda_1 \cdot \sigma_{\tau} (\blambda_{-1}) & \geq & 
\lambda_2 \cdot \sigma_{\tau} (\blambda_{-2}) \;\; \geq \;\;
\ldots
\;\; \geq \;\;
\lambda_n \cdot \sigma_{\tau} (\blambda_{-n}) .
\ea
\eeq
Indeed, without loss of generality, let us justify the first inequality:
$$
\ba{rcl}
\lambda_1 \cdot \sigma_{\tau} (\blambda_{-1}) & \geq & 
\lambda_2 \cdot \sigma_{\tau} (\blambda_{-2}) 
\ea
$$
Recall that 
$$
\ba{rcl}
\sigma_{\tau}(\aa) & \Def & \sum\limits_{1 \leq i_1 < \ldots < i_{\tau} \leq n - 1}
a_{i_1} \cdot \ldots \cdot a_{i_{\tau}}, \qquad \aa \in \R^{n - 1}.
\ea
$$
Hence,
$$
\ba{rcl}
\lambda_1 \cdot \sigma_{\tau} (\blambda_{-1})
& = & \lambda_1 \cdot 
\Bigl[ \,
\sum\limits_{2 \leq i_1 < \ldots < i_{\tau} \leq n}
\lambda_{i_1} \cdot \ldots \cdot \lambda_{i_{\tau}} \, \Bigr] \\
\\
& = &
\lambda_1 \cdot \lambda_2 \cdot
\Bigl[ 
\,
\sum\limits_{ 3 \leq i_2 < \ldots < i_{\tau} \leq n}
\lambda_{i_2} \cdot \ldots \cdot \lambda_{i_{\tau}} 
\, \Bigr]
\;\; + \;\; 
\lambda_1 \cdot 
\Bigl[
\,
\sum\limits_{ 3 \leq i_1 < \ldots < i_{\tau} \leq n}
\lambda_{i_1} \cdot \ldots \cdot \lambda_{i_{\tau}}
\,
\Bigr] \\
\\
& \geq & 
\lambda_1 \cdot \lambda_2 \cdot
\Bigl[ 
\,
\sum\limits_{ 3 \leq i_2 < \ldots < i_{\tau} \leq n}
\lambda_{i_2} \cdot \ldots \cdot \lambda_{i_{\tau}} 
\, \Bigr]
\;\; + \;\; 
\lambda_2 \cdot 
\Bigl[
\,
\sum\limits_{ 3 \leq i_1 < \ldots < i_{\tau} \leq n}
\lambda_{i_1} \cdot \ldots \cdot \lambda_{i_{\tau}}
\,
\Bigr] \\
\\
& = &
\lambda_2 \cdot \Bigl[   
\, 
\sum\limits_{
\substack{ 1 \leq i_1 \leq \ldots \leq \lambda_{i_{\tau}} \leq n   \\
				  2 \notin \{ i_1, \ldots, i_{\tau} \} }}   
			  \lambda_{i_1} \cdot \ldots \cdot \lambda_{i_{\tau}}
\,\Bigr]
\;\; = \;\; \lambda_2 \cdot \sigma_{\tau} (\blambda_{-2}).
\ea
$$
Therefore, we have established \eqref{th4_2_bound} and obtain:
$$
\ba{rcl}
\lambda_n \cdot \sigma_{\tau} (\blambda_{-n}) \mat{I}
\;\;\; \preceq \;\;\;
\mat{B}^{1/2} \mat{P}_{\tau} \mat{B}^{1/2}
& \preceq &
\lambda_1 \cdot \sigma_{\tau} (\blambda_{-1}) \mat{I},
\ea
$$
which proves the required bound. \qed

\subsection{Proof of \cref{TheoremVolumeSampling}}

Let $\mat{B} = \mat{Q} \DiagonalMatrix(\Vector{\lambda}) \mat{Q}\Transpose$
be a spectral decomposition of~$\mat{B}$,
where $\Vector{\lambda} = \Vector{\lambda}(\mat{B})$
and $\mat{Q} \in \RealMatrices{n}{n}$ is an orthogonal matrix.
Formula~(3.5) in~\cite{rodomanov2020randomized} states that
\[
  \Expectation_{S \DistributedAs \VolumeSampling_{\tau + 1}(\mat{B})} [
    \mat{I}_S (\mat{B}_{S \times S})^{-1} \mat{I}_S\Transpose
  ]
  =
  \frac{1}{\sigma_{\tau + 1}(\Vector{\lambda})}
  \mat{Q}
  \DiagonalMatrix\bigl(
    \sigma_{\tau}(\Vector{\lambda}_{-1}),
    \ldots,
    \sigma_{\tau}(\Vector{\lambda}_{-n})
  \bigr)
  \mat{Q}\Transpose.
\]
(Their $\sigma_{\tau}$ is $\sigma_{\tau + 1}$ in our notation.)
But, according to \cref{LemmaSpec},
\[
  \mat{Q}
  \DiagonalMatrix\bigl(
    \sigma_{\tau}(\Vector{\lambda}_{-1}),
    \ldots,
    \sigma_{\tau}(\Vector{\lambda}_{-n})
  \bigr)
  \mat{Q}\Transpose
  =
  \mat{P}_{\tau},
\]
and the claim follows.
\qed

\subsection{Proof of Lemma~\ref{LemmaXiProperties}}

Clearly, when $\tau = 0$, we have $\xi_{0}(\blambda) \equiv 1$ .
For $\tau = n - 1$,
inequalities \eqref{PApproxQuality} are in fact identities, and
$\xi_{n - 1}(\blambda)  \overset{\eqref{PN_1}}{=}
\frac{\lambda_n}{\lambda_1}.$

Now, let us prove that $\xi_{\tau}(\blambda)$
monotonically decrease with $\tau$.
For that, we fix some $1 \leq \tau < \rho \leq n - 1$
and justify 
$$
\ba{rcl}
\xi_{\tau}(\blambda) & \Def & \frac{\sigma_{\tau}(\blambda_{-1})}{\sigma_{\tau}(\blambda_{-n})}
\;\; \geq \;\;
\xi_{\rho}(\blambda) 
\;\; \Def \;\; \frac{\sigma_{\rho}(\blambda_{-1})}{\sigma_{\rho}(\blambda_{-n})},
\ea
$$
which is equivalent to
\beq \label{SigmaMonotonicity}
\ba{rcl}
\sigma_{\tau}(\blambda_{-1}) \cdot \sigma_{\rho}(\blambda_{-n})
& \geq & 
\sigma_{\rho}(\blambda_{-1}) \cdot \sigma_{\tau}(\blambda_{-n}).
\ea
\eeq
Then, \eqref{SigmaMonotonicity} can be rewritten as 
\beq \label{Monotone2}
\ba{rcl}
(A + \lambda_n B) \cdot (C + \lambda_1 D)
& \geq & (A + \lambda_1 B ) \cdot (C + \lambda_n D), 
\ea
\eeq
where
$$
\ba{rcl}
A & \Def & 
\sum\limits_{2 \leq i_1 < \ldots < i_{\tau} \leq n - 1}
\lambda_{i_1} \cdot \ldots \cdot \lambda_{i_{\tau}}, \\
\\
B & \Def & 
\sum\limits_{2 \leq i_1 < \ldots < i_{\tau - 1} \leq n - 1}
\lambda_{i_1} \cdot \ldots \cdot \lambda_{i_{\tau - 1}}, \\
\\
C & \Def & 
\sum\limits_{2 \leq i_1 < \ldots < i_{\rho} \leq n - 1}
\lambda_{i_1} \cdot \ldots \cdot \lambda_{i_{\rho}}, \\
\\
D & \Def & 
\sum\limits_{2 \leq i_1 < \ldots < i_{\rho - 1} \leq n - 1}
\lambda_{i_1} \cdot \ldots \cdot \lambda_{i_{\rho - 1}}.
\ea
$$
By opening up the brackets in \eqref{Monotone2} and eliminating
the common terms, we obtain
$$
\ba{rcl}
\lambda_1 (AD - BC) & \geq & \lambda_n (AD - BC).
\ea
$$
Thus, to justify \eqref{SigmaMonotonicity}, it is enough to prove that $AD - BC \geq 0$,
which is the following inequality:
$$
\ba{rcl}
 L_1 & \Def & \sum\limits_{
	\substack{
	2 \leq i_1 < \ldots < i_{\tau} \leq n - 1 \\
	2 \leq j_1 < \ldots < j_{\rho - 1} \leq n - 1}}
\lambda_{i_1} \cdot \ldots \cdot \lambda_{i_{\tau}}
\cdot \lambda_{j_1} \cdot \ldots \cdot \lambda_{j_{\rho -1}} \\
& \geq &
\sum\limits_{
	\substack{
		2 \leq i_1 < \ldots < i_{\tau - 1} \leq n - 1 \\
		2 \leq j_1 < \ldots < j_{\rho} \leq n - 1}}
\lambda_{i_1} \cdot \ldots \cdot \lambda_{i_{\tau - 1}}
\cdot \lambda_{j_1} \cdot \ldots \cdot \lambda_{j_{\rho}}
\;\; \Def \;\; L_2,
\ea
$$
which is true since $\rho > \tau$.
Indeed, we can rewrite the left hand side sum as
$$
\ba{rcl}
L_1 & = & \sum\limits_{k = \tau + 1}^{n - 1} \lambda_k \cdot
\biggl[ \, 
\sum\limits_{\substack{
2 \leq i_1 < \ldots < i_{\tau - 1} < k \\
2 \leq j_1 < \ldots < j_{\rho - 1} \leq n - 1
}}
\lambda_{i_1} \cdot \ldots \cdot \lambda_{i_{\tau - 1}}
\cdot \lambda_{j_1} \cdot \ldots \cdot \lambda_{j_{\rho - 1}} 
\, \Bigr] \\
& \geq & 
\sum\limits_{k = \rho + 1}^{n - 1} \lambda_k \cdot
\biggl[ \, 
\sum\limits_{\substack{
		2 \leq i_1 < \ldots < i_{\tau - 1} < k \\
		2 \leq j_1 < \ldots < j_{\rho - 1} \leq n - 1
}}
\lambda_{i_1} \cdot \ldots \cdot \lambda_{i_{\tau - 1}}
\cdot \lambda_{j_1} \cdot \ldots \cdot \lambda_{j_{\rho - 1}} 
\, \Bigr] \\
& \geq &
\sum\limits_{k = \rho + 1}^{n - 1} \lambda_k \cdot
\biggl[ \, 
\sum\limits_{\substack{
		2 \leq i_1 < \ldots < i_{\tau - 1} < n - 1 \\
		2 \leq j_1 < \ldots < j_{\rho - 1} \leq k
}}
\lambda_{i_1} \cdot \ldots \cdot \lambda_{i_{\tau - 1}}
\cdot \lambda_{j_1} \cdot \ldots \cdot \lambda_{j_{\rho - 1}} 
\, \Bigr] \\
& = & L_2,
\ea
$$
where we used that $\rho > \tau$.

Finally, to prove the limit, let us divide the right hand side of
$$
\ba{rcl}
\xi_{\tau}(\blambda) & = & \frac{\sigma_{\tau}(\blambda_{-1})}{\sigma_{\tau}(\blambda_{-n})}
\;\; \leq \;\;
\frac{ \sum\limits_{2 \leq i_1 < \ldots < i_{\tau} \leq n - 1} \lambda_{i_1} \cdot \ldots
	\cdot \lambda_{i_{\tau}}  }{  
\lambda_1 \cdot \ldots \cdot \lambda_{\tau}}
\ea
$$
by the biggest element from the sum,
which is $\lambda_2 \cdot \ldots \cdot \lambda_{\tau + 1}$.
Thus, we get
$$
\ba{rcl}
\xi_{\tau}(\blambda) & \leq & 
\frac{1 + E}{(\lambda_1 / \lambda_{\tau + 1})},
\ea
$$
where $E$ is a finite sum of numbers that are smaller than or equal to $1$. Hence, 
$\xi_{\tau}(\blambda) \to 0$ when $\frac{\lambda_1}{\lambda_{\tau + 1}} \to \infty$. \qed

\subsection{Proof of Proposition~\ref{PropositionChebyshev}}
\label{SubsectionAppendixChebyshev}

Let us use an upper bound on $\gamma$ from \eqref{GammaPDef},
which is the \textit{uniform} polynomial approximation for the whole interval
$[\lambda_n, \lambda_1]$:
\beq \label{GammaChebyshev}
\ba{rcl}
\!\!\!
\gamma(p_{\tau}) & \leq &
\max\limits_{s \in [\lambda_n, \lambda_1]} | sp_{\tau}(s) - 1 |.
\ea
\eeq

Then, we use
\beq \label{QChoice}
\ba{rcl}
Q_{\tau}(s) & \Def &
T_{\tau + 1}\bigl( \frac{\lambda_1 + \lambda_n - 2s}{\lambda_1 - \lambda_n} \bigr)
\cdot  T_{\tau + 1}\bigl( \frac{\lambda_1 + \lambda_n}{\lambda_1 - \lambda_n} \bigr)^{-1},
\ea
\eeq
where $T_{\tau + 1}(\cdot)$ is the standard Chebyshev polynomial of the first kind of degree $\tau + 1$.
Namely, we can define them recursively:
$$
\ba{rcl}
T_{0}(x) & \Def & 1,
\qquad 
T_{1}(x) \;\; \Def \;\; x, 
\qquad 
T_{k + 1}(x) \;\; \Def \;\; 2 x \cdot T_{k}(x) - T_{k - 1}(x), \qquad k \geq 1.
\ea
$$

Note that $Q_{\tau}(0) = 1$, thus the polynomial $1 - Q_{\tau}(s)$ 
is divisible by $s$. Then, we take 
$$
\ba{rcl}
p_{\tau}(s) & := & \frac{1 - Q_{\tau}(s)}{s},
\ea
$$
which is the polynomial of degree $\tau$. 
	This choice ensures that
	$$
	\ba{rcl}
	\gamma & \overset{\eqref{GammaChebyshev},\eqref{QChoice}}{\leq} &
	\max\limits_{x \in [-1, 1]} | T_{\tau + 1}(x) |
	\cdot
	T_{\tau + 1}\bigl( \frac{\lambda_1 + \lambda_n}{\lambda_1 - \lambda_n} \bigr)^{-1} \\[10pt]
	& = &
	T_{\tau + 1}\bigl( \frac{\lambda_1 + \lambda_n}{\lambda_1 - \lambda_n} \bigr)^{-1}
	\; \leq \;
	2 \Bigl( \frac{\sqrt{\lambda_1} - \sqrt{\lambda_n}}
	{\sqrt{\lambda_1} + \sqrt{\lambda_n}}  \Bigr)^{\tau + 1},
	\ea
	$$
	where the last inequality is the classical bound for the Chebyshev polynomials
	(see, e.g. Section 16.4 in \cite{vishnoi2013lx}).
	Thus, the condition number 
	$$
	\ba{rcl}
	\frac{\beta}{\alpha} & = & \frac{1 + \gamma}{1 - \gamma}
	\ea
	$$
	decreases 
	exponentially with $\tau$.
	\qed
\subsection{Proof of Theorem~\ref{TheoremKrylovGM}}

Let us fix an arbitrary $\mat{P} = \mat{P}^{\top} \succ 0$
such that $\mat{P} = p_\tau(\mat{B})$, for some
polynomial $p_{\tau} \in \R[s]$ and $\deg(p_{\tau}) = \tau$.
We take $\beta := \beta(\mat{P})$ and $\alpha := \alpha(\mat{P})$ 
(from the definition \eqref{PBound})
and denote 
$$
\ba{rcl}
\mat{\bar{P}}  & := & \frac{1}{\beta L} \mat{P}.
\ea
$$

Let us consider an arbitrary iteration $k \geq 0$, and denote
the following step
$$
\ba{rcl}
\mat{T} & := & \xx_k - \mat{\bar{P}} \nabla f(\xx_k).
\ea
$$
Recall also that 
$$
\ba{rcl}
\xx_{k + 1} & := & \xx_k - \mat{P}_{\aa_k} \nabla f(\xx_k).
\ea
$$
By the optimality of $\aa_k$
as the projection of $\mat{B}^{-1} \nabla f(\xx_k)$ onto the Krylov subspace, we have:
\beq \label{a_k_proj_bound}
\ba{cl}
& \la \nabla f(\xx_k), \xx_{k + 1} - \xx_k \ra + \frac{L}{2} \| \xx_{k + 1} - \xx_k \|_{\mat{B}}^2
\;\; \equiv \;\;
\frac{L}{2} \| \xx_{k + 1} - \xx_k + \frac{1}{L}\mat{B}^{-1} \nabla f(\xx_k) \|_{\mat{B}}^2
- \frac{1}{2L} \| \nabla f(\xx_k) \|_{\mat{B}^{-1}}^2 \\
\\
& \; \leq \;
\frac{L}{2} \| \mat{T} - \xx_k + \frac{1}{L}\mat{B}^{-1} \nabla f(\xx_k) \|_{\mat{B}}^2
- \frac{1}{2L} \| \nabla f(\xx_k) \|_{\mat{B}^{-1}}^2
\;\; \equiv \;\;
\la \nabla f(\xx_k), \mat{T}  - \xx_k \ra + \frac{L}{2} \| \mat{T} - \xx_k \|_{\mat{B}}^2.
\ea
\eeq
Hence, we obtain, for any $\yy \in \R^n$:
$$
\ba{rcl}
f(\xx_{k + 1}) & \overset{\eqref{FuncGlobalBounds}}{\leq} &
f(\xx_k) + \la \nabla f(\xx_k), \xx_{k + 1} - \xx_k \ra + \frac{L}{2}\| \xx_{k + 1} - \xx_k\|_{\mat{B}}^2 \\
\\
& \overset{\eqref{a_k_proj_bound}}{\leq} & 
f(\xx_k) + \la \nabla f(\xx_k), \mat{T} - \xx_k \ra + \frac{L}{2}\| \mat{T} - \xx_k\|_{\mat{B}}^2 \\
\\
& \leq & 
f(\xx_k) + \la \nabla f(\xx_k), \mat{T} - \xx_k \ra + \frac{\beta L}{2}\| \mat{T} - \xx_k\|_{\mat{P}^{-1}}^2 \\
\\
& \leq & 
f(\xx_k) + \la \nabla f(\xx_k), \yy - \xx_k \ra + \frac{\beta L}{2} \|\yy - \xx_k\|_{\mat{P}^{-1}}^2.
\ea
$$
where we used that $\mat{T}$ is the minimizer of the last upper bound in $\yy$.
Thus, using convexity, we get, for any $\yy \in \R^n$:
\beq \label{KrylovUpperBound}
\ba{rcl}
f(\xx_{k + 1}) & \leq & f(\yy) + \frac{\beta L}{2} \|\yy - \xx_k \|_{P^{-1}}^2
\;\; \leq \;\;
f(\yy) +  \frac{\beta}{\alpha} \cdot \frac{L}{2} \|\yy - \xx_k \|_{\mat{B}}^2.
\ea
\eeq
In particular, substituting $\yy := \xx_k$ we justify that the method is \textit{monotone}:
$f(\xx_{k + 1}) \leq f(\xx_k), \forall k \geq 0$. Therefore,
$$
\ba{rcl}
\xx_k & \in & \mathcal{F}_{0} \;\; \Def \;\; \Bigl\{  
\, \xx \in \R^n \; : \; f(\xx) \leq f(\xx_0),
\Bigr\}
\ea
$$
and we assume that the initial level set $\mathcal{F}_0$ is \textit{bounded}, denoting
$$
\ba{rcl}
D_0 & \Def & \sup\limits_{\xx\in \mathcal{F}_0} \|\xx - \xx^{\star} \|_{\mat{B}} \;\; < \;\; +\infty.
\ea
$$
Substituting $\yy := \gamma_k \xx^{\star} + (1 - \gamma_k) \xx_k$, $\gamma_k \in [0, 1]$
into \eqref{KrylovUpperBound}, we obtain
\beq \label{Th5_1MainBound}
\ba{rcl}
f(\xx_{k + 1}) & \leq & \gamma_k f^{\star} + (1 - \gamma_k) f(\xx_k)
+ \gamma_k^2 \frac{\beta}{\alpha} \cdot \frac{L}{2}\|\xx^{\star} - \xx_k \|_{\mat{B}}^2 \\
\\
& \leq & 
\gamma_k f^{\star} + (1 - \gamma_k) f(\xx_k)
+ \gamma_k^2 \frac{\beta}{\alpha} \cdot \frac{L}{2}D_0^2.
\ea
\eeq
Substituting $\gamma_k := \frac{2}{k + 1}$ and using the standard technique (see the proof of Theorem~\ref{TheoremGM}), we establish the global rate for the convex case:
$$
\ba{rcl}
f(\xx_k) - f^{\star} & \leq & \cO\Bigl( \,  \frac{\beta}{\alpha} \cdot \frac{L D_0^2}{k}  \, \Bigr).
\ea
$$
For strongly convex functions ($\mu > 0$), we continue as
$$
\ba{rcl}
f(\xx_{k + 1}) & \overset{\eqref{Th5_1MainBound}, \eqref{FuncGlobalBounds}}{\leq} &
\gamma_k f^{\star} + (1 - \gamma_k) f(\xx_k)
+ \gamma_k^2 \frac{\beta L }{\alpha \mu} \cdot \bigl( f(\xx_k) - f^{\star}  \bigr),
\ea
$$
and choosing $\gamma_k := \frac{\alpha \mu}{2 \beta L}$ we establish the exponential rate. \qed
\section{Adaptive Search}
\label{SectionAdaptive}

In this section, we briefly present adaptive versions of \cref{alg:GM,alg:FGM}
which do not require the knowledge of the constant $M = \beta L$ and
can automatically ``tune'' it in iterations yet preserving the original
worst-case efficiency estimates.
This is achieved by using a standard ``backtracking line search'' which
can be found, e.g., in~\cite{nesterov2013gradient}.

In what follows, for any $\Vector{x}, \Vector{y} \in \EffectiveDomain \psi$,
$M > 0$ and $\mat{P} = \mat{P}\Transpose \succ 0$,
we define the following predicate:
\[
  \QuadraticGrowth_{M, \mat{P}}(\Vector{x}, \Vector{y}) \colon
  \quad
  f(\Vector{y})
  \leq
  f(\Vector{x})
  +
  \InnerProduct{\Gradient f(\Vector{x})}{\Vector{y} - \Vector{x}}
  +
  \frac{M}{2} \RelativeNorm{\Vector{y} - \Vector{x}}{\mat{P}^{-1}}^2.
\]
According to our assumptions~\eqref{FuncGlobalBounds} and~\eqref{PBound},
we know that this predicate is surely satisfied for any pair of points once
$M \geq \beta L$.

The adaptive version of \cref{alg:GM} is presented in \cref{alg:AdaptiveGM}.
This method starts with a certain initial guess~$\tilde{M}_0$ for the
constant~$\beta L$ and then, at every iteration, repeatedly increases the
current guess in two times until the predicate becomes satisfied.
This process is guaranteed to terminate (when $M_k$ becomes bigger or equal to
$\beta L$, or even sooner).
After that, we accept the new point $x_{k + 1}$ and choose a new ``optimistic''
guess of the constant~$M$ for the next iteration by halving the value of $M_k$
that we have accepted at the current iteration.

\begin{algorithm}
  \caption{Adaptive Preconditioned GM}
  \label{alg:AdaptiveGM}
  \begin{algorithmic}
    \STATE \textbf{Input:}
      $\Vector{x}_0 \in \EffectiveDomain \psi$,
      $\mat{P} = \mat{P}^{\top} \succ 0$,
      $\tilde{M}_0 > 0$.
    \FOR{$k = 0, 1, \dots$}
      \STATE
        Find smallest integer $i_k \geq 0$ such that
        \[
          \Vector{x}_{k + 1}
          =
          \GradientStep_{M_k, \mat{P}}\bigl(
            \Vector{x}_k, \Gradient f(\Vector{x}_k)
          \bigr),
          \quad
          M_k = 2^{i_k} \tilde{M}_k
        \]
        satisfies the predicate
        $\QuadraticGrowth_{M_k, \mat{P}}(\Vector{x}_k, \Vector{x}_{k + 1})$.
      \STATE
        Set $\tilde{M}_{k + 1} = M_k / 2$.
    \ENDFOR
  \end{algorithmic}
\end{algorithm}

We assume that the preconditioner~$\mat{P}$ is sufficiently simple so that
we can efficiently check the predicate
$\QuadraticGrowth_{M_k, \mat{P}}(\Vector{x}_k, \Vector{x}_{k + 1})$.
For example, if $\psi = 0$, then
$
  \Vector{x}_{k + 1}
  =
  \Vector{x}_k - \frac{1}{M_k} \mat{P} \Gradient f(\Vector{x}_k)
$
and
$
  M_k \RelativeNorm{\Vector{x}_{k + 1} - \Vector{x}_k}{\mat{P}^{-1}}^2
  =
  \InnerProduct{\Gradient f(\Vector{x}_k)}{\Vector{x}_k - \Vector{x}_{k + 1}}
$
can be efficiently computed.

For \cref{alg:AdaptiveGM}, we can prove exactly the same rates as
in \cref{TheoremGM} (up to absolute constants) provided that
\begin{equation}
  \label{eq:RequirementOnInitialGuess}
  \tilde{M}_0 \leq \beta L.
\end{equation}
The proof is essentially the same as in \cref{sec:ProofGM} with only two
minor differences: 1) inequality~\eqref{Th1OneStep} is now guaranteed by our
predicate; 2) instead of using $M = \beta L$ in~\eqref{Th1PreTelescoped},
we should use the bound $M_k \leq 2 \beta L$ which follows
from~\eqref{eq:RequirementOnInitialGuess} and the fact that any value
of $M \geq \beta L$ is always acceptable in the line search.
Using a classical argument from~\cite{nesterov2013gradient}, it is not difficult
to show that, on average, \cref{alg:AdaptiveGM} makes only $\sim 2$ gradient
steps at each iteration.

In contrast to an upper estimate of the constant $\beta L$, an initial guess
satisfying~\eqref{eq:RequirementOnInitialGuess} can be easily generated.
One simple recipe is to make a trial step
$
  \Vector{x}_1'
  =
  \GradientStep_{M_0', \mat{P}}\bigl(
    \Vector{x}_0, \Gradient f(\Vector{x}_0)
  \bigr)
$
for an \emph{arbitrarily chosen} $M_0' > 0$ and then compute
\[
  \tilde{M}_0
  =
  \frac{
    f(\Vector{x}_1')
    -
    f(\Vector{x}_0)
    -
    \InnerProduct{\Gradient f(\Vector{x}_0)}{\Vector{x}_1' - \Vector{x}_0}
  }{
    \frac{1}{2} \RelativeNorm{\Vector{x}_1' - \Vector{x}_0}{\mat{P}^{-1}}^2
  }.
\]
Alternatively, we can find a suitable $\tilde{M}_0$ be choosing
an arbitrary $M_0' > 0$ and then repeatedly halving it until the predicate
$\QuadraticGrowth(\Vector{x}_0, \Vector{x}_1'(M))$ stops being satisfied for
$
  \Vector{x}_1'(M)
  =
  \GradientStep_{M, \mat{P}}(\Vector{x}_0, \Gradient f(\Vector{x}_0))
$.
This auxiliary procedure either terminates in a logarithmic number of steps,
in which case we get a suitable $\tilde{M}_0$, or, otherwise, we quickly find
an approximate solution of our problem.

Similar technique can be applied for the Fast Gradient Method.
Specifically, let us introduce an auxiliary procedure shown in
\cref{alg:FastGradientStep} for computing one iteration of \cref{alg:FGM}
for a given value of~$M$.
Then, the adaptive FGM method can be constructed as shown in \cref{alg:AdaptiveFGM}.
As in the basic method, we can show that the rates from \cref{TheoremFGM}
still remain valid (up to absolute constants) for \cref{alg:AdaptiveFGM},
provided that $\tilde{M}_0$ satisfies~\eqref{eq:RequirementOnInitialGuess}.
For generating the initial guess $\tilde{M}_0$, we can use exactly the same
techniques as before.

\begin{algorithm}
  \caption{
    $
      (\Vector{x}_+, \Vector{v}_+, A_+; \Vector{y})
      =
      \FastGradientStep_{M, \rho, \mat{P}}(\Vector{x}, \Vector{v}, A)
    $
  }
  \label{alg:FastGradientStep}
  \begin{algorithmic}
    \REQUIRE
      $M > 0$;
      $\rho \geq 0$;
      $\mat{P} = \mat{P}\Transpose \succ 0$;
      $\Vector{x}, \Vector{v} \in \EffectiveDomain \psi$;
      $A > 0$.
    \STATE
      Find $a_+$ from eq.\
      $\frac{M a_+^2}{A + a_+} = 1 + \rho (A + a_+)$.
    \STATE
      Set $A_+ = A + a_+$,
      $H = \frac{1 + \rho A_+}{a_+}$,
      $\theta = \frac{a_+}{A_+}$,
      $\omega = \frac{\rho}{H}$,
      $\gamma = \frac{\omega (1 - \theta)}{1 - \omega \theta}$.
    \STATE
      Set
      $\hat{\Vector{v}} = (1 - \gamma) \Vector{v} + \gamma \Vector{x}$,
      $\Vector{y} = (1 - \theta) \Vector{x} + \theta \hat{\Vector{v}}$.
    \STATE
      Compute
      $
        \Vector{v}_+
        =
        \GradientStep_{H, \mat{P}}\bigl(
          \hat{\Vector{v}}, \nabla f(\Vector{y})
        \bigr)
      $.
    \STATE
      Set
      $
        \Vector{x}_+
        =
        (1 - \theta) \Vector{x} + \theta \Vector{v}_+
      $.
    \RETURN $(\Vector{x}_+, \Vector{v}_+, A_+; \Vector{y})$.
  \end{algorithmic}
\end{algorithm}

\begin{algorithm}
  \caption{Adaptive Preconditioned FGM}
  \label{alg:AdaptiveFGM}
  \begin{algorithmic}
    \STATE \textbf{Input:}
      $\Vector{{x}}_0 \in \EffectiveDomain \psi$,
      $\mat{P} = \mat{P}\Transpose \succ 0$, $\tilde{M}_0 > 0$.
    \STATE
      Set $\Vector{v}_0 = \Vector{x}_0$, $A_0 = 0$.
    \FOR{$k = 0, 1, \dots$}
      \STATE
        Find smallest integer $i_k \geq 0$ such that
        \[
          (\Vector{x}_{k + 1}, \Vector{v}_{k + 1}, A_{k + 1}; \Vector{y}_k)
          =
          \FastGradientStep_{M_k, \mat{P}}\bigl(
            \Vector{x}_k, \Vector{v}_k, A_k
          \bigr),
          \qquad
          M_k = 2^{i_k} \tilde{M}_k
        \]
        satisfies the predicate
        $\QuadraticGrowth_{M_k, \mat{P}}(\Vector{y}_k, \Vector{x}_{k + 1})$.
      \STATE
        Set $\tilde{M}_{k + 1} = M_k / 2$.
    \ENDFOR
  \end{algorithmic}
\end{algorithm}

\end{document}